\documentclass{amsproc}

\usepackage[utf8]{inputenc}
\usepackage{amsmath}
\usepackage{amsthm}
\usepackage{amsfonts}
\usepackage{amssymb}
\usepackage{nicefrac}
\usepackage{xcolor}
\usepackage{graphicx}
\usepackage{caption}
\usepackage{subcaption}
\usepackage{hyperref}
\usepackage{tikz}
\usepackage{forest}
\usepackage{upgreek}
\usepackage{blkarray}
\usepackage{graphics}
\usepackage{relsize}
\usepackage{graphicx}
\usepackage{nicematrix}
\usepackage{algorithm}
\usepackage{algpseudocode}
\usepackage{pifont}

\newcommand{\centered}[1]{\begin{tabular}{c} #1 \end{tabular}}

\newcommand{\RR}{\mathbb{R}}
\newcommand{\PP}{\mathbb{P}}
\newcommand{\CC}{\mathbb{C}}
\newcommand{\rk}{\mathrm{rk}}
\newcommand{\Gr}{\mathrm{Gr}}

\algnewcommand{\OfflineInputs}[1]{%
  \State \textbf{Offline Input}
  \Statex \hspace*{\algorithmicindent}\parbox[t]{.8\linewidth}{\raggedright #1}
}

\algnewcommand{\OnlineInputs}[1]{%
  \State \textbf{Online Input}
  \Statex \hspace*{\algorithmicindent}\parbox[t]{.8\linewidth}{\raggedright #1}
}

\algnewcommand{\OfflinePre}[1]{%
  \State \textbf{Offline Precomputations}
  \Statex \hspace*{\algorithmicindent}\parbox[t]{.8\linewidth}{\raggedright #1}
}

\algnewcommand{\OnlineComp}[1]{%
  \State \textbf{Online Computations}
  \Statex \hspace*{\algorithmicindent}\parbox[t]{.8\linewidth}{\raggedright #1}
}

\algnewcommand{\OnlineOut}[1]{%
  \State \textbf{Online Output}
  \Statex \hspace*{\algorithmicindent}\parbox[t]{.8\linewidth}{\raggedright #1}
}

\algnewcommand{\Prelim}[1]{%
  \State \textbf{Precomputations for a Random Instance}
  \Statex \hspace*{\algorithmicindent}\parbox[t]{.8\linewidth}{\raggedright #1}
}

\algnewcommand{\BasisB}[1]{%
  \State \textbf{Sample a Monomial Basis $\mathcal{B}$}
  \Statex \hspace*{\algorithmicindent}\parbox[t]{.8\linewidth}{\raggedright #1}
}

\algnewcommand{\ElimE}[1]{%
  \State \textbf{Compute a Syzygy-Reduced Elimination Template  $\mathcal{E}$}
  \Statex \hspace*{\algorithmicindent}\parbox[t]{.8\linewidth}{\raggedright #1}
}

\algnewcommand{\TempOfflineOut}[1]{%
  \State \textbf{Possible Candidate for Output}
  \Statex \hspace*{\algorithmicindent}\parbox[t]{.8\linewidth}{\raggedright #1}
}

\algnewcommand{\OfflineOut}[1]{%
  \State \textbf{Offline Output}
  \Statex \hspace*{\algorithmicindent}\parbox[t]{.8\linewidth}{\raggedright #1}
}

\newtheorem{theorem}{Theorem}[section]
\newtheorem{lemma}[theorem]{Lemma}

\theoremstyle{definition}
\newtheorem{definition}[theorem]{Definition}
\newtheorem{example}[theorem]{Example}

\theoremstyle{remark}
\newtheorem{remark}[theorem]{Remark}

\numberwithin{equation}{section}

\begin{document}

\title{Snapshot of Algebraic Vision}


%    author one information
\author{Joe Kileel}
\address{Department of Mathematics and Oden Insitute, University of Texas at Austin}
\curraddr{}
\email{\href{mailto:jkileel@math.utexas.edu}{jkileel@math.utexas.edu}}
\thanks{J.K. is supported in part by NSF awards DMS-2309782 and IIS-2312746, and  start-up grants from the Department of Mathematics and Oden Institute at UT Austin.}

%    author two information
%\author{Kathl\'en Kohn \\ \begin{center} \begin{small} \textit{In honor of Bernd Sturmfels' 60th birthday} \end{small} \end{center}}
\author{Kathl\'en Kohn}
\address{Institutionen f\"or Matematik, KTH Royal Institute of Technology}
\curraddr{}
\email{\href{mailto:kathlen@kth.se}{kathlen@kth.se}}
\thanks{K.K. is supported in part by the Wallenberg AI, Autonomous Systems and Software Program (WASP) funded by the Knut and Alice Wallenberg Foundation.}

%    The 2020 edition of the Mathematics Subject Classification is
%    the current definitive version.

%\date{{In honor of Bernd Sturmfels' 60th birthday}}

\begin{abstract}
 In this survey article, we present interactions between  algebraic geometry and computer vision, which have recently come under the header of algebraic vision. 
The subject has given new insights in multiple view geometry and its application to 3D scene reconstruction
and carried a host of novel  problems and ideas back into algebraic geometry.
\end{abstract}

\dedicatory{In honor of Bernd Sturmfels' 60th birthday}

\maketitle

\textit{Computer vision} is the research field that studies how computers can gain understanding from 2D images and videos, similar to human cognitive abilities.
Typical computer vision tasks include the automatic recognition of objects in images, the detection of events in videos, and the reconstruction of 3D scenes from many given 2D images.
A general overview of computer vision is presented in  textbook form in~\cite{szeliski2010computer}.  
The subject is a pillar in the AI revolution.

\textit{Algebraic vision} is the symbiosis of computer vision and  algebraic geometry. 
Motivated by Chris Aholt's Ph.D. thesis titled \emph{Polynomials in Multiview Geometry} \cite{aholtThesis} and earlier works,
the term ``algebraic vision" was coined during a particular lunch held at a Seattle office of Google in early spring 2014, attended by
Sameer Agarwal, Chris Aholt, Joe Kileel, Hon-Leung Lee, Max Lieblich, Bernd Sturmfels, and Rekha Thomas.
The intent was to encourage interactions between the applied algebraic geometry community and the 3D reconstruction community in computer vision.  
A short discussion of algebraic vision can be found in the review \cite{breiding2021nonlinear} on nonlinear algebra and its applications.

Historically, computer vision  made  substantial use of projective geometry and computational algebra in parts of its foundations.
Specifically \textit{multiple view geometry}, as described in the textbook \cite{HZ} of Hartley and Zisserman, is modeled on projective three-space and two-space and group-equivariant (multi-)linear transformations between these.  
Similar algebraic treatments of the subject are the textbooks \cite{ma2004invitation} and \cite{maybank2012theory}.
Previously, this connection was not well-appreciated by the computational algebra geometry community.
However, in the last decade, algebro-geometric papers and workshops on 3D reconstruction have been appearing, leading to novel results in multiple view geometry while motivating developments in applied algebraic geometry.

The present article provides a  survey  of algebraic vision.  
No previous knowledge of computer vision is assumed, and the prerequisites for computational algebraic geometry are kept mostly to the level of undergraduate texts  \cite{cox2013ideals}. 
Due to space limitations, the article makes no attempt to be comprehensive in any way, but instead it focuses narrowly on the role of projective varieties and systems of polynomial equations in 3D vision. 
An outline of the sections is as follows:
\begin{itemize}
\item In Section~\ref{sec:computervisionproblems}, we introduce the problem of 3D scene reconstruction from unkown cameras and its algebro-geometric nature.
\item In Section~\ref{sec:cameramodels}, we discuss a variety of usual models for cameras. 
\item In Section~\ref{sec:multiviewvarieties}, we study \textit{multiview varieties} which
 characterize feasible images of points under fixed cameras. 
Their defining equations play a key role in 3D reconstruction algorithms, and their Euclidean distance degrees measure the intrinsic complexity of noisy triangulation (i.e., the task of recovering the 3D coordinates of a point observed by known cameras).
\item In Section~\ref{sec:cameraTuples}, we consider the space of all cameras. 
We explain how tuples of cameras -- up to changes of world coordinates -- can be encoded via \textit{multifocal tensors} \cite{HZ}.
\item In Section~\ref{sec:3Dalgorithms}, we overview the most popular algorithmic pipeline to solve 3D scene reconstruction, highlighting \textit{minimal problems} that are the algebro-geometric heart of the pipeline.

\item In Section \ref{ssec:eliminationTemplates}, we describe polynomial solvers for minimal problems, focusing on Gröbner basis methods using \textit{elimination templates} and homotopy continuation.  Those method applies to zero-dimensional parameterized polynomial systems in general. 

\item In Section~\ref{sec:discriminants}, we discuss algebro-geometric approaches to understand degenerate world scenes and image data, where uniqueness of reconstruction breaks down and algorithms can encounter difficulty.

\end{itemize}

After reading Sections \ref{sec:computervisionproblems} and \ref{sec:cameramodels}, the other sections are essentially independent; only Section \ref{ssec:eliminationTemplates} builds on Section \ref{sec:3Dalgorithms}.
We provide specific pointers to earlier sections in case of partial dependencies.

Some important topics in algebraic vision that are omitted include group synchronization (e.g., \cite{rosen2019se,lerman2021robust}), uses of polynomial optimization  (e.g., \cite{kahl2007globally,cifuentes2022local,aholt2012qcqp, yang2022certifiably, ROB-077}), and approaches based on differential invariants (e.g., \cite{calabi1998differential, boutin2004structure}).   Readers may consult \cite{ozyecsil2017survey} for a survey that covers numerical and large-scale optimization aspects in 3D reconstruction.

\medskip

 \textbf{Acknowledgements.} 
We thank Sameer Agarwal,  Paul Breiding, Luca Carlone, Tim Duff, Hongyi Fan, Fredrik Kahl, Anton Leykin, Tomas Pajdla, Jean Ponce, Kristian Ranestad, Felix Rydell, Elima Shehu, Rekha Thomas, Matthew Trager and Uli Walther for their comments on earlier versions of the manuscript.

\section{Computer vision  through the algebraic lens} \label{sec:computervisionproblems}
One of the main challenges in computer vision is the \emph{structure-from-motion} (SfM) problem:
given many 2D images, the task is to reconstruct the 3D scene and also the positions of the cameras that took the pictures.
This has many applications such as 3D mapping from images taken by drones \cite{remondino2012uav}, to localize and navigate autonomous cars and robots in a 3D world \cite{grisetti2010tutorial}, or
in the movie industry to reconstruct 3D backgrounds \cite{Kitagawa:Mocap:book}, for photo tourism \cite{agarwal2011building}, and for combining real and virtual worlds \cite{Dobbert:Matchmoving:book}.

\begin{figure}[htb]
    \centering
     \begin{subfigure}[b]{0.21\textwidth}
         \centering
         \includegraphics[width=\textwidth]{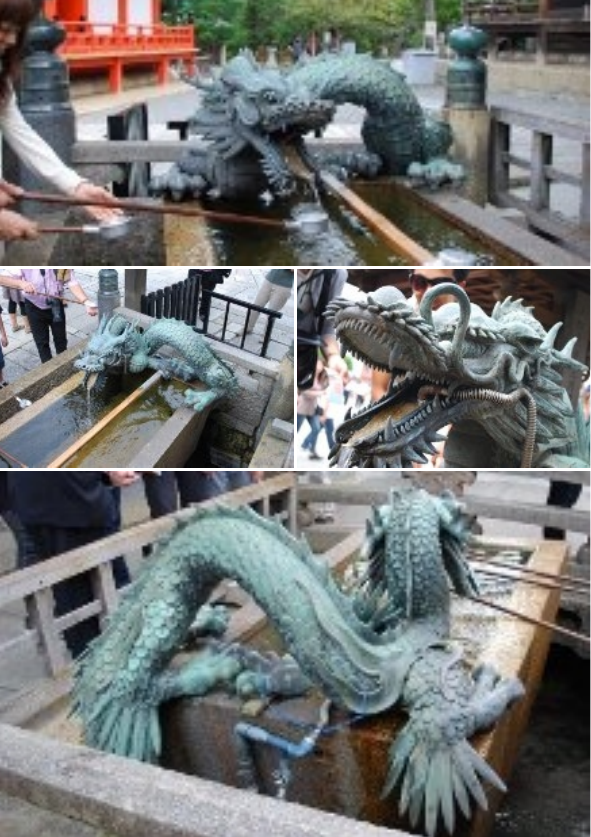}
         \caption{Input images}
         \label{fig:input}
     \end{subfigure}
          \begin{subfigure}[b]{0.23\textwidth}
         \centering
         \includegraphics[width=0.92\textwidth]{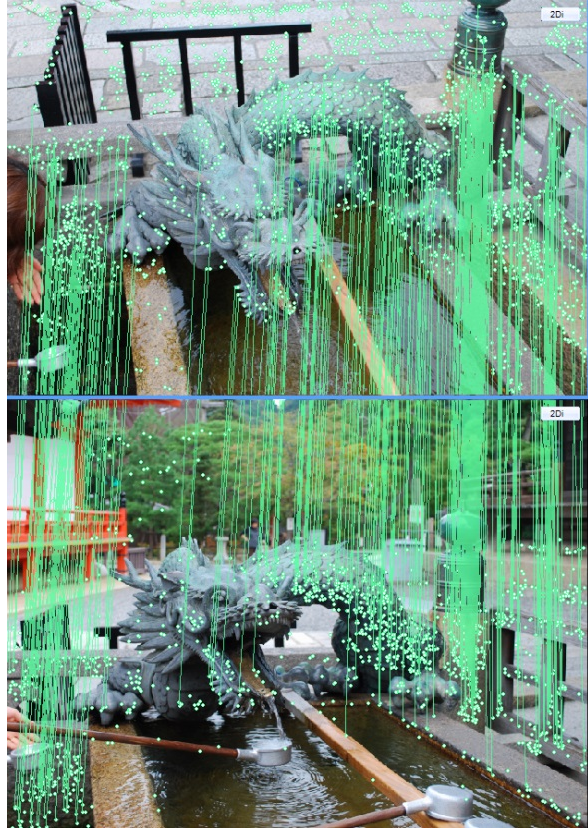}
         \caption{Image matching}
         \label{fig:matching}
     \end{subfigure}
          \begin{subfigure}[b]{0.25\textwidth}
         \centering
         \includegraphics[width=\textwidth]{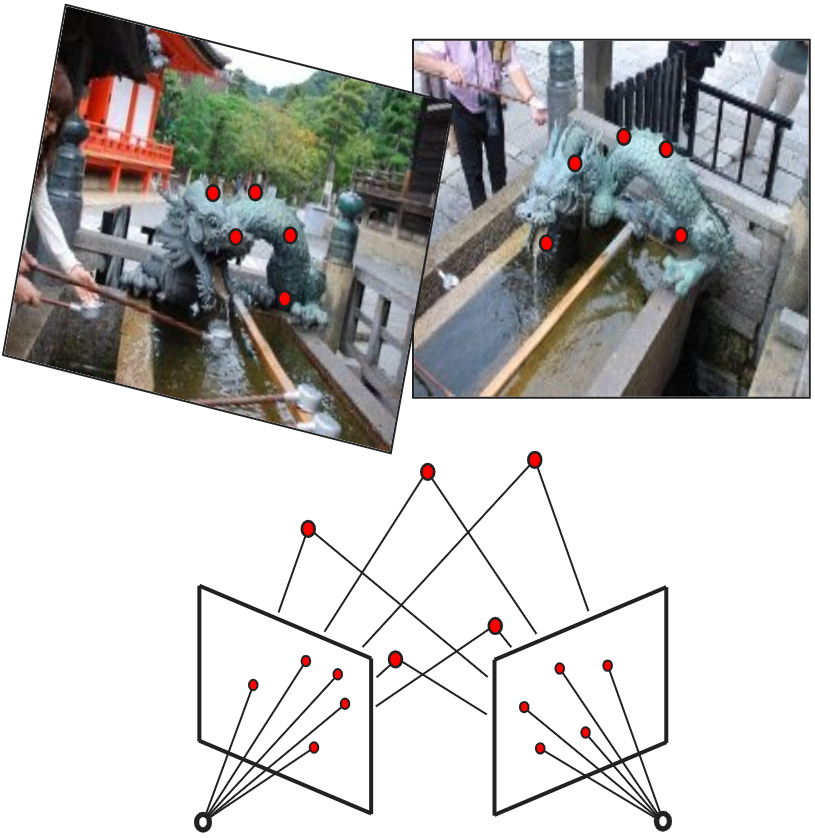}
         \captionsetup{justification=raggedleft
}
         \caption{Reconstruct cam- eras and 3D points}
         \label{fig:algebraStep}
     \end{subfigure}
          \begin{subfigure}[b]{0.28\textwidth}
         \centering
         \includegraphics[width=\textwidth]{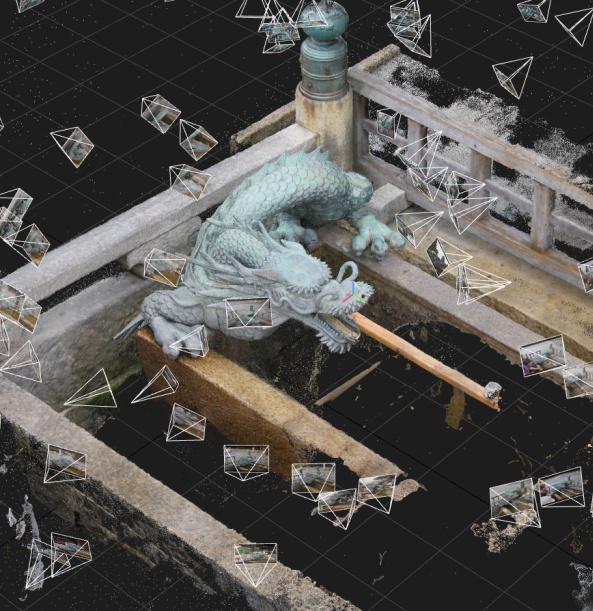}
         \caption{Output}
         \label{fig:output}
     \end{subfigure}
    \caption{3D reconstruction pipeline (courtesy of Tomas Pajdla).}
    \label{fig:pipeline}
\end{figure}

The {structure-from-motion} problem is typically solved using the 3D reconstruction pipeline. 
We will now sketch a highly simplified version of that pipeline, illustrated in Figure~\ref{fig:pipeline}.
We provide more details in Section \ref{ssec:pipeline}.
Given a set of 2D images,
the first step in the pipeline is to take a few of the given images and identify geometric features, such as points or lines,  that they have in common.
In Figure~\ref{fig:matching}, a detection algorithm has been used that only identifies points.
In the second step of the pipeline, we forget the original images and only keep the geometric features we have identified.
We reconstruct the 3D coordinates of those features and also the camera poses, that is, the locations and orientations of the cameras. 
In Figure~\ref{fig:algebraStep}, five common points were identified on two images, so we aim to reconstruct the five points in 3-space and the two cameras.
Finally, we repeat this process several times until we have recovered all cameras and also enough geometric features to approximate the 3D scene.

As the second step of the pipeline forgets the pictures and only works with algebro-geometric features, such as points or lines, the reconstruction problem becomes purely algebraic. 
More specifically, we aim to compute a fiber of the \emph{joint camera map:}
\begin{align}
    \label{eq:jointCameraMap}
    \Phi: \mathcal{X} \times \mathcal{C}_m \dashrightarrow \mathcal{Y},
\end{align}
that maps an arrangement $X \in \mathcal{X}$ of 3D features
and a tuple $(C_1, \ldots, C_m) \in \mathcal{C}_m$ of cameras 
to the $m$ 2D images of $X$ taken by the cameras.
For instance in Figure~\ref{fig:algebraStep}, the joint camera map becomes
\begin{align} \label{eq:jointCameraMapExampleAffine}
    \Phi: \left(  \RR^3 \right)^5 \times \mathcal{C}_2 \dashrightarrow
    \left(  \RR^2 \right)^5 \times \left(  \RR^2 \right)^5.
\end{align}

\begin{figure}[htb]
    \centering
    \includegraphics[width=0.48\textwidth,height=5.4cm]{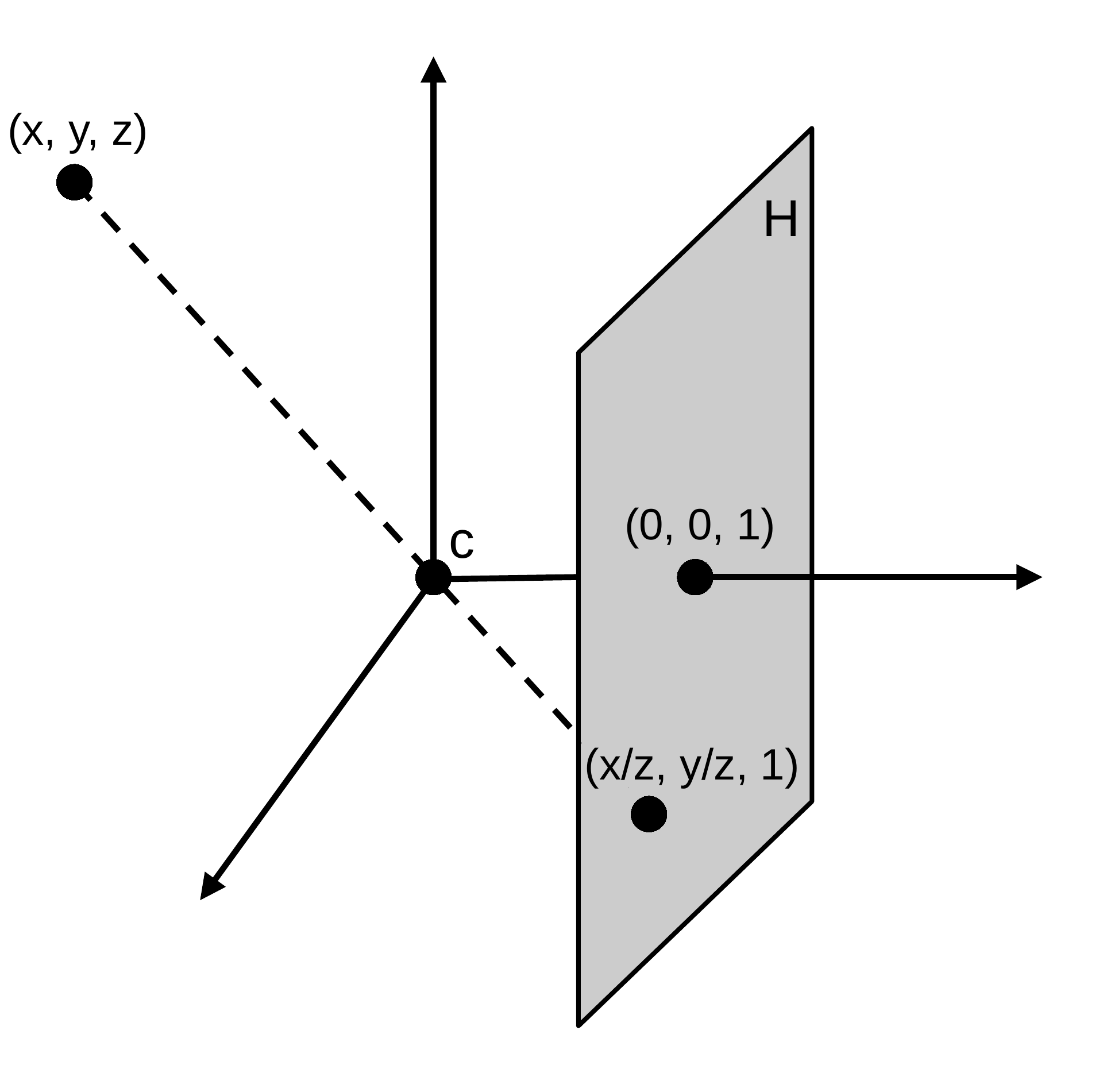}
    \caption{A pinhole camera  in standard position is centered at  $c=(0,0,0)$ and maps world points $(x,y,z)$ to image points $( \tfrac{x}{z},  \tfrac{y}{z}, 1)$ on the image plane $H$.}
    \label{fig:pinhole}
\end{figure}

A full specification of the joint camera map requires a choice of camera model.
The simplest  model is a \emph{pinhole camera}; see Figure~\ref{fig:pinhole}.
Such a camera simply takes a picture of a point in space by projecting it onto a plane.
A pinhole camera in standard position is typically assumed to be centered at the origin such that its image plane is
$H = \{ (x,y,z) \in \RR^2 \mid z=1 \}$.
In these coordinates, the pinhole camera is the map
\begin{align*}
    \RR^3 \dashrightarrow H, \quad
    (x,y,z) \longmapsto ( \tfrac{x}{z},  \tfrac{y}{z}, 1).
\end{align*}

Often homogeneous coordinates are used to model cameras.
This means that each point in the image plane is identified with the light ray passing through the point and the origin. 
In homogeneous coordinates, the standard pinhole camera in Figure \ref{fig:pinhole} becomes
\begin{align*}
    \PP_{\RR}^3 \dashrightarrow \PP_{\RR}^2, \quad 
    [x:y:z:w] \longmapsto [x:y:z].
\end{align*}
This map is defined everywhere except at the camera center $[0:0:0:1]$, i.e. the origin in the affine chart where $w = 1$.

The projective geometry approach in modeling cameras is thoroughly explained in the textbook \cite{HZ}.
That book laid many foundations and conventions used in modern computer vision and offers a great entry point for the algebraic community into the field of computer vision.
The main focus of the book \cite{HZ} is multiview geometry, where a 3D object is viewed by several cameras, such as in Figure \ref{fig:pipeline}. 
In that setting, we cannot assume that all cameras are in standard position as described above. 
Instead, a pinhole camera is more generally given by a $3 \times 4$ matrix $A$ of rank three.
The corresponding camera map 
\begin{align*}
    \PP_{\RR}^3 \dashrightarrow \PP_{\RR}^2, \quad 
    X \longmapsto AX.
\end{align*}
is defined everywhere except at the camera center that is given by the kernel of~$A$.
The standard camera in Figure~\ref{fig:pinhole} corresponds to the matrix $\left[ \begin{smallmatrix}  1 & 0 & 0 & 0 \\ 0 & 1 & 0 & 0 \\ 0 & 0 & 1 & 0 \end{smallmatrix} \right]$.

Hence, when using  pinhole cameras and homogeneous coordinates, the \emph{camera variety} $\mathcal{C}_m$ in \eqref{eq:jointCameraMap}
that describes all $m$-tuples of such cameras is $$ \mathcal{C}_m=
(\PP\,\mathrm{Mat}^{3 \times 4}_3 )^m,$$ where $\mathrm{Mat}^{3 \times 4}_3 \subset \mathbb{R}^{3 \times 4}$ denotes the set of $3 \times 4$ matrices of rank three. 
For instance, the joint camera map in \eqref{eq:jointCameraMapExampleAffine} becomes
\begin{align*}
    \Phi: \left(  \PP_\RR^3 \right)^5 \times (\PP\,\mathrm{Mat}^{3 \times 4}_3 )^2 &\,\dashrightarrow
    \left(  \PP_\RR^2 \right)^5 \times \left(  \PP_\RR^2 \right)^5, \\
    (X_1, \ldots, X_5, A_1, A_2) &\longmapsto (A_1 X_1, \ldots, A_1 X_5, A_2 X_1, \ldots, A_2 X_5).
\end{align*}

In the next section, we review common camera models and highlight algebraic  vision articles studying camera geometry.
In the remaining sections, our focus returns to the joint camera map in \eqref{eq:jointCameraMap}:
We will see that many computer vision problems can be formulated using the joint camera map~--
such as understanding the image of a shape in space or reconstructing a 3D shape from several images~-- and are thus natural to study through the algebraic lens.
The recent paper \cite{agarwal2022atlas} gives a similar such unifying algebro-geometric framework for computer vision problems.

\section{Camera models} \label{sec:cameramodels}

\textbf{Calibrated cameras.}
The camera model described in the previous section is  known as the \emph{projective / uncalibrated pinhole camera}.
The \emph{calibrated pinhole camera} model
assumes that every camera is obtained from the standard pinhole camera in Figure~\ref{fig:pinhole} by translation and rotation.
This means that every camera matrix $A$ is of the form $[R \mid t]$ where $R \in \mathrm{SO}(3)$ is the relative rotation from the standard pinhole camera to the camera with matrix $A$ and the relative translation can be read off from the vector $t \in \mathbb{R}^3$:
the camera center $c$, which is the origin in Figure~\ref{fig:pinhole}, is now $c = - R^{\top} t$ (note that the vector $(c,1)^\top \in \mathbb{R}^4$ spans the kernel of the camera matrix $[R \mid t]$).
In particular, every calibrated pinhole camera has 6 degrees of freedom (3 for $R$ and 3 for $t$), whereas a projective pinhole camera has 11 degrees of freedom.

Calibrated pinhole cameras are a commonly used model in applications, corresponding to the case when the internal parameters of the cameras are known (such as from meta data stored inside the image file). 
There is also a variety of \emph{partly calibrated pinhole cameras}, e.g. a camera with unknown focal length, that have less strict structural assumptions on the $3 \times 4$ camera matrices than the fully calibrated model described above. 
Partly calibrated pinhole cameras are modeled as $K [R \mid t]$ where $K$ is a $3 \times 3$ upper triangular \emph{calibration matrix} whose entries are partially known \cite[Chapter 6]{HZ}. 

\begin{figure}
    \centering
    \includegraphics[width = 0.3\textwidth]{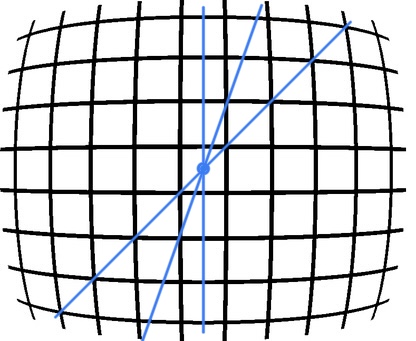}
    \caption{A grid affected by radial distortion. Radial lines and distortion center are blue.}
    \label{fig:radial}
\end{figure}

\textbf{Distortion.}
In practice, cameras are not as ideal as in the calibrated model. 
As seen in Figure~\ref{fig:pinhole}, the pinhole cameras described so far assume that the world point, the camera center, and the image point are collinear.
This assumption does not hold for real-life camera lenses, because they are affected by various kinds of distortion. 
The main factor of deviation from the idealistic pinhole camera model is typically radial distortion; see Figure~\ref{fig:radial}. 

Often, calibrated cameras are a sufficient approximation of real-life cameras.
However, sometimes the impact of radial distortion is too big, e.g., for fisheye cameras. 
One approach to address radial distortion is to make the camera model more complicated by adding distortion parameters that have to be estimated during 3D reconstruction
(see \cite[Chapter 7.4]{HZ} for an overview and \cite{kileel2018distortion} for an algebraic treatment of distortion varieties).
Another approach is to simplify the camera model by not estimating the radial distortion at all: 
Once the center of radial distortion on a given image is determined, we know for every 3D point onto which line through the distortion center it gets mapped by the camera (see Figure \ref{fig:radial}), although we do not know its exact point in the image. 
Thus, the camera simplifies to a map $\PP^3 \dashrightarrow \PP^1$ that sends 3D points to radial lines through the distortion center in the image. 
This is known as the \emph{1D radial camera model} \cite{larsson2020calibration,thirthala2012radial}.

\textbf{Rolling shutter and pushbroom cameras.}
Many real-world cameras do not function like pinhole cameras that were 
described so far. For instance, smartphones, that produce a massive and 
yet growing amount of picture data, use \emph{rolling shutters} (see e.g., \cite{dai2016rolling}). When taking a 
picture, these cameras scan across the scene, capturing the resulting 
image row by row. For static cameras taking pictures of still scenes, 
there is no difference between rolling shutters and standard pinhole 
cameras that capture an entire image at the same instant (also known as 
global shutter). However, for moving cameras or objects, in particular 
videos, global shutter pinhole cameras are not a good approximation of 
rolling shutters. For example, when a rolling shutter camera moves with 
constant speed along a line, the picture of any other line in 
three-space would typically be a conic in the image plane, whereas a 
pinhole camera maps lines to lines.

A related camera
model is the \textit{pushbroom camera}, which is common in 
satellites \cite[Chapter 6.4.1]{HZ}. Such a camera considers only a plane in 
the 3D world that passes through the camera center and then projects 
that plane onto a single line. Hence, a non-moving pushbroom camera 
produces only an image line, not a plane. Two-dimensional images arise 
by moving the camera.

Rolling shutter and pushbroom cameras have received little attention 
from the algebraic community so far.  
In particular, the theory described in Sections~ \ref{sec:multiviewvarieties}--\ref{sec:discriminants} is largely undeveloped for such cameras.

\textbf{Congruences.}
Ponce, Sturmfels, and Trager \cite{ponce2017congruences} propose a general framework to study camera geometry that includes all camera models without distortion discussed above. Their notion of a camera is an abstraction from the classical way of modeling a camera by a map $\PP^3 \dashrightarrow \PP^2 $. Instead of identifying a concrete image plane with coordinates, they define a \emph{rational geometric camera} to be a map $\PP^3 \dashrightarrow \mathrm{Gr}(1, \PP^3)$ from three-space into the Grassmannian of lines. The fiber over any image point under a classical camera $\PP^3 \dashrightarrow \PP^2 $ without distortion is a line, namely a \emph{viewing ray} of the camera. The associated rational geometric camera maps each world point directly to its corresponding viewing ray. For instance, a geometric pinhole camera assigns to every world point distinct from the camera center the unique line spanned by the world point and the camera center. 

The image of a rational geometric camera is a \emph{congruence}, i.e. a two-dimensional family of lines. To obtain an actual photographic camera, one composes the geometric camera map with a map from its image congruence to $\PP^2$ that models the image formation process via mapping viewing rays to image coordinates. 
These physical realizations of rational geometric cameras are the focus of Trager, Sturmfels, Canny, Hebert, and Ponce in \cite{trager2017general}.

In \cite{ponce2017congruences}, the authors suggest an even more general camera definition by saying that a \emph{geometric camera} is an arbitrary congruence in $\mathrm{Gr}(1, \PP^3)$. The \emph{order} of a congruence is the number of lines in the congruence that pass through a general point in $\PP^3$. Hence, congruences of order one are exactly those that are the images of rational geometric cameras $\PP^3 \dashrightarrow \mathrm{Gr}(1, \PP^3)$. Congruences of higher order are geometric cameras such that a world point typically appears several times on each image taken by the camera. This happens for instance for many moving rolling shutter cameras, or for \textit{catadioptric cameras} which are cameras that use curved mirrors.
The idea of modeling cameras by congruences had been expressed by computer vision researchers before~\cite{pajdla2002stereo,ponce2009camera,sturm2005multi}.

\textbf{Chirality.}
In real life, a 3D point has to lie in front of a camera in order to be seen by it.
Imposing this semialgebraic constraint on a camera model is referred to as \emph{chirality}.
The theoretical investigation of chirality was initiated by Hartley \cite{hartley1998chirality}; see also \cite[Chapter 21]{HZ}.
Further studies of this concept were for instance undertaken by Laveau and Faugeras \cite{laveau1996oriented}, Werner and Pajdla \cite{werner2001cheirality,werner2001oriented,werner2003combinatorial,werner2003constraint}, and Agarwal, Pryhuber, Sinn, and Thomas \cite{agarwal2020chiral,pryhuber2020existence}.

\section{Multiview varieties} \label{sec:multiviewvarieties}
Given $m$ fixed pinhole cameras $C = (C_1, \ldots, C_m)$, where $C_i: \PP_\RR^3 \dashrightarrow \PP_\RR^2$, we can specialize the joint camera map $\Phi$ from \eqref{eq:jointCameraMap} to
\begin{align}
\label{eq:jointImageMap}
    \Phi_{C}: \PP_\RR^3 \dashrightarrow (\PP_\RR^2)^m,
\end{align}
sending a 3D point to its $m$ images under the cameras.
The \emph{joint image} of the cameras, first introduced by Triggs \cite{triggs1995matching}, is $\Phi_C(\PP_\RR^3)$.
The first algebraic study of the joint image was done by Heyden and Åström \cite{heyden1997algebraic}. They called it the \emph{natural descriptor} and observed that it is not Zariski closed. 
The Zariski closure $\mathcal{M}_C := \overline{\Phi_C(\PP_\RR^3)}$ of the joint image is called the \emph{joint image variety} \cite{trager2015joint} or \emph{multiview variety} \cite{aholt2013hilbert}.
It coincides with the Euclidean closure of the joint image
\cite[Theorem 4]{agarwal2020chiral}.

\subsection{Multiview constraints} \label{sec:multiviewconstraints}
Many computer vision works have studied systems of polynomials that vanish on the joint image, often referred to as \emph{multiview  constraints}; 
see e.g.
\cite{faugeras1995geometry,heyden1997algebraic,faugeras2001geometry,ma2004invitation,aholt2013hilbert,trager2015joint,agarwal2019ideals}.
Their main motivation is that the joint image completely determines the cameras $C_1, \ldots, C_m$, up to changes of coordinates in $\PP^3$.
So one way of thinking about the 3D reconstruction pipeline in Figure~\ref{fig:pipeline} is to first obtain point correspondences between the given images (i.e., (noisy) samples on the joint image) and then estimate sufficiently many multiview constraints from the correspondences so as to be able to reconstruct the cameras, then the 3D points.

We now describe the multiview constraints studied by Heyden and Åström \cite{heyden1997algebraic}.
Let $A_1, \ldots, A_m$ be $3 \times 4$ matrices of rank three that define the cameras, i.e., $C_i: X \mapsto A_iX$.
We write $x_i = (x_{i1}:x_{i2}:x_{i3})$ for the homogeneous coordinates on the $i$-th image $\PP_\RR^2$, 
and form the $3m \times (m+4)$ matrix
\begin{align*}
    M_A = M_A(x_1, \ldots, m) := \begin{bmatrix} A_1 & x_1 & 0 & \cdots & 0 \\
     A_2 & 0 & x_2 & \cdots & 0 \\
     \vdots & \vdots & \vdots & \ddots & \vdots \\
     A_m & 0 & 0 & \cdots & x_m\end{bmatrix}.
\end{align*}
A tuple of image points $(x_1, \ldots, x_m) \in (\PP_\RR^2)^m$ lies in the joint image if and only if there is some 3D point $X$ and non-zero scalars $\lambda_i$ such that $A_i X = \lambda_i x_i$ for all $i =1, \ldots, m$. This implies that $[X,-\lambda_1, \ldots, -\lambda_m]^\top$ is in the kernel of $M_A$, so the matrix $M_A$ must be rank deficient.
In fact, the multiview variety $\mathcal{M}_C$ is cut out by the maximal minors of $M_A$:
\begin{align}
\label{eq:mLinearities}
    \mathcal{M}_C = \{ (x_1, \ldots, x_m) \in (\PP_\RR^2)^m \mid \rk\, M_A(x_1, \ldots, x_m) < m+4 \}.
\end{align}
These maximal minors are multilinear, i.e. they are linear in every triplet $x_i$ of variables.

\begin{example} \label{ex:fundamendalMatrix}
For $m=2$ cameras, $M_A$ is a $6 \times 6$ matrix.
Its determinant defines the multiview hypersurface $\mathcal{M}_C$ in $\PP_\RR^2 \times \PP_\RR^2$.
This equation is a bilinear form in $x_1, x_2$ represented by a matrix $F \in \RR^{3 \times 3}$, i.e. $\det(M_A) = x_2^{\top} F x_1$.  Here $F$ is called the \emph{fundamental matrix} for projective pinhole cameras, and the \emph{essential matrix} for calibrated pinhole cameras. 
Longuet-Higgins' seminal work \cite{longuet1981computer} introduced the essential matrix to the computer vision community along with a first 3D reconstruction algorithm.
We will discuss essential matrices, fundamental matrices, and their generalizations to more than two cameras in Section \ref{sec:cameraTuples}.
\end{example}

For every subset $J \subseteq \{ 1, \ldots, m \}$ of cardinality $k$, we consider the $3k \times (k+4)$ submatrix $M_{A,J}$ of $M_A$ whose rows and columns correspond to the cameras $C_j$ for $j \in J$.
As the maximal minors of $M_{A,J}$ are multilinear constraints in the $k$ triplets of variables $\{ x_j \mid j \in J \}$, 
we refer to the set of all maximal minors of all $M_{A,J}$ with $|J|=k$ as the \emph{$k$-linearities}.
As seen in \eqref{eq:mLinearities}, the multiview variety $\mathcal{M}_C$ of $m$ cameras is cut out by the $m$-linearities.
Heyden and Åström \cite{heyden1997algebraic} showed that $\mathcal{M}_C$ is in fact cut out by the $2$-linearities if the $m$ camera centers are not coplanar.

An ideal-theoretic treatment of the multiview variety was first done by Aholt, Sturmfels, and Thomas \cite{aholt2013hilbert},  
and later continued by Agarwal, Pryhuber, and Thomas \cite{agarwal2019ideals}.
They study the vanishing ideal of $\mathcal{M}_C$, called the \emph{multiview ideal}.

\begin{theorem}
\begin{enumerate}
    \item If the $m$ camera centers are pairwise distinct, the multiview ideal is generated by the $2$- and $3$-linearities \cite[Theorem 3.7]{agarwal2019ideals}.
    \item If all $4 \times 4$ minors of the $4 \times 3m$ matrix $[ A_1^\top \; A_2^\top \; \cdots \; A_m^\top ]$ are non-zero, then the $2$-, $3$-, and $4$-linearities form a  universal Gröbner basis of the multiview ideal \cite[Theorem 2.1]{aholt2013hilbert}.
\end{enumerate}
\end{theorem}

A nice overview of the joint image and its Zariski closure (equivalently, Euclidean closure) is presented in \emph{the joint image handbook} by Trager, Hebert, and Ponce \cite{trager2015joint}.
They provide an explicit description of the actual joint image (instead of its closure) as a constructible set, and discuss  sets of $k$-linearities that cut out the multiview variety as well as  sets of $k$-linearities that are not sufficient to define $\mathcal{M}_C$ but uniquely determine the cameras $C_1, \ldots, C_m$, up to coordinates changes in~$\PP_{\mathbb{R}}^3$.

\subsection{Euclidean distance degree and triangulation} \label{ssec:edd}
\emph{Triangulation} refers to the problem of reconstructing the 3D coordinates of a point from 2D images when the cameras are known.
For  pinhole cameras $C=(C_1, \ldots, C_m)$,  this means to compute the fiber of the specialized joint camera map $\Phi_C$ in \eqref{eq:jointImageMap}.
 In theory, triangulation is trivial to solve for $m \geq 2$ cameras,  since the fiber under each camera $C_i: \PP_\RR^3 \dashrightarrow \PP_\RR^2$ is a line and so the desired 3D point is the intersection of these $m$ lines in $\PP_\RR^3$.
 However, in practice, the image measurements are noisy, which means that the given image points $(x_1, \ldots, x_m) \in (\PP_\RR^2)^m$ do not lie on the multiview variety $\mathcal{M}_C$.
 Thus, to solve triangulation, one first needs to \textit{denoise} the given data, for instance by finding a nearby point on the multiview variety $\mathcal{M}_C$ to the given data point $x = (x_1, \ldots, x_m) \in (\PP_\RR^2)^m$.
 
 In practical applications, the given image points are finite, which means that after a choice of affine chart $\RR^2$ of each $\PP_\RR^2$ we may assume that $x \in \RR^{2m}$.
 In this setting, we want to find the closest point $x$ on the 
 restriction of $\mathcal{M}_C$ to the chosen affine chart, that is, on the \emph{affine multiview variety}: 
 $$\mathcal{M}_C^\alpha = \mathcal{M}_C \cap \RR^{2m}.$$
 Typically, distance is measured via the squared Euclidean distance on $\RR^{2m}$.
 Regarding $\mathcal{M}_C^\alpha$ as a complex variety in $\CC^{2m}$, for almost all data points $x \in \CC^{2m}$, this optimization problem has a constant number of complex critical points \cite{draisma2016euclidean}. 
 That number is known as the \emph{Euclidean distance degree} (\emph{ED degree} for short) of the variety $\mathcal{M}_C^\alpha$.
ED degree measures the algebraic complexity of the triangulation problem, in the sense that it provides an upper bound for the number of real critical points for generic data points $x \in \RR^{2m}$. 
 
 Hartley and Sturm \cite{hartley1997triangulation} showed that the ED degree of the affine multiview variety $\mathcal{M}_{(C_1,C_2)}^\alpha$ for two cameras is $6$.
 Stew{\'e}nius and Nist{\'e}r  computed the ED degree of $\mathcal{M}_C^\alpha$ for $m \leq 7$ cameras and observed that it grows cubically \cite{hartley2007optimal}. 
 The problem of determining the ED degree of the affine multiview variety actually motivated the general introduction of Euclidean distance degrees of algebraic varieties by Draisma, Horobe{\c{t}}, Ottaviani, Sturmfels, and Thomas \cite[Example 3.3]{draisma2016euclidean}, 
 as well as follow-up work on orthogonally invariant matrix sets \cite{drusvyatskiy2015counting,drusvyatskiy2017euclidean}.
 In \cite[Conjecture 3.4]{draisma2016euclidean}, a precise cubic formula for the ED degree of  $\mathcal{M}_C^\alpha$ for an arbitrary number of cameras was conjectured.
 Harris and Lowengrub \cite{harris2018chern} provided an upper bound for the ED degree, showing that it is indeed bounded by some cubic polynomial in the number of cameras, 
 before \cite[Conjecture 3.4]{draisma2016euclidean} was finally proven by Maxim, Rodriguez, and Wang:
 \begin{theorem}[{\cite{maxim2020euclidean}}]
When the $m \geq 2$ cameras are in general position, the Euclidean distance degree of the affine multiview variety $\mathcal{M}_C^\alpha$ is $$\frac{9}{2}m^3 - \frac{21}{2}m^2 + 8m - 4.$$
 \end{theorem}
\noindent Informally, the theorem means triangulation with a large number of cameras may be quite difficult: the number of critical points grows cubically with $m$.

\subsection{Extensions}
There are several extensions of multiview varieties which incorporate additional constraints or apply to different models. 
We discuss some recent extensions below.

\textbf{Chirality.}
The discussion of multiview varieties and triangulation above ignores the constraint of chirality, i.e., that the 3D point must lie in front of each camera in order to be seen by it.
Agarwal, Pryhuber, Sinn, and Thomas \cite{agarwal2020chiral} provide semialgebraic descriptions (using polynomial equalities and inequalities) of both the \emph{chiral domain}, that is the set of points in $\PP_\RR^3$ that lie in front of each of the $m$ pinhole cameras $C_i$, and the \emph{chiral joint image}, i.e. the image of the chiral domain under the specialized joint camera map $\Phi_C$ in \eqref{eq:jointImageMap}.
The chiral joint image is a subset of Triggs' joint image  $\Phi_C(\PP_\RR^3)$ \cite{triggs1995matching} and describes the true image of the 3D world as seen by the cameras $C=(C_1, \ldots, C_m)$.
Hence, the semialgebraic description of the chiral joint image from \cite{agarwal2020chiral} is a refinement of the multiview constraints discussed above.

\textbf{Congruences.}
Multiview varieties can be naturally generalized to other camera models (besides pinhole cameras as described above).
Recall that a rational geometric camera is a map
\begin{align*}
    \PP^3 \dashrightarrow \Sigma \subset \Gr(1,\PP^3)
\end{align*}
that maps a world point $X$ to a line passing through $X$.
The Zariski closure of the image of the map is a congruence $\Sigma$, i.e. a surface in the Grassmannian $\Gr(1,\PP^3)$ of lines. 
The map is defined everywhere except at its \emph{focal locus}, that is, the set of  points $X \in \PP^3$ that belong to more than one line on $\Sigma$.

Given $m$ rational geometric cameras $C = (C_1, \ldots, C_m)$,  $C_i: \PP^3 \dashrightarrow \Sigma_i \subset \Gr(1,\PP^3)$, 
the \emph{multi-image variety} $\mathcal{M}_C$ is the Zariski closure of the image of
\begin{align*}
    \Phi_C: \PP^3 &\,\dashrightarrow \Sigma_1 \times \ldots \times \Sigma_m \subset \Gr(1, \PP^3)^m, \quad
    X \longmapsto \left( C_1(X), \ldots, C_m(X) \right).
\end{align*}
Ponce, Sturmfels, and Trager \cite{ponce2017congruences} showed that the multi-image variety often has an easy description in terms of the \emph{concurrent lines variety} $V_m \subset \Gr(1, \PP^3)^m$ that consists of all $m$-tuples of lines that meet in a common point.

\begin{theorem}[{\cite[Theorem 5.1]{ponce2017congruences}}]
If the $m$ focal loci of the cameras $C_1, \ldots, C_m$ are pairwise disjoint, then 
$\mathcal{M}_C = V_m \cap  (\Sigma_1 \times \ldots \times \Sigma_m)$.
\end{theorem}

The authors also describe a minimal set of generators and a reduced Gröbner bases for the prime ideal of the concurrent lines variety $V_m$ in \cite[Theorem 3.1]{ponce2017congruences}.
Moreover, using the Pl\"ucker embedding $\Gr(1, \PP^3) \subset \PP^5$, they conjectured a formula for the multidegree of $V_m$ in $(\PP^5)^m$, which was proven by Escobar and Knutson \cite[Theorem 2.4]{escobar2017multidegree}. The latter theorem also provides the multidegree of the multi-image variety $\mathcal{M}_C$ in $(\PP^5)^m$ under the assumption that the $m$ focal loci of the cameras $C_1, \ldots, C_m$ are pairwise disjoint.

\textbf{3D data.}  
Another way to generalize multiview varieties is by not only considering the image of a single 3D point under the given cameras, but rather instead considering the image of more complex 3D objects.
The \emph{line multiview variety}, i.e., the Zariski closure of the image of the map $\Gr(1, \PP^3) \dashrightarrow (\Gr(1, \PP^2))^m$ that sends a 3D line to its images under $m$ pinhole cameras, is described in \cite{breiding2022line,faugeras1995geometry}.
Multiview varieties for more than one point are studied in  \cite{joswig2016rigid}. 
For instance, the authors there consider the set of point pairs in $\RR^3 \times \mathbb{R}^3$ that are a Euclidean distance of $1$ apart from each other.
Let $\mathcal{X} \subset \PP_\RR^3 \times \PP_\RR^3$ be the Zariski closure of said set.
Given $m$ pinhole cameras $C = (C_1, \ldots, C_m)$, 
the authors restrict the joint camera map from \eqref{eq:jointCameraMap} to:
\begin{align}\label{eq:rigid-variation}
    \mathcal{X} &\,\dashrightarrow (\PP_\RR^2)^m \times (\PP_\RR^2)^m, \quad
    (X_1,X_2) \longmapsto \left( \Phi_C(X_1), \Phi_C(X_2) \right),
\end{align}
and define the \emph{rigid multiview variety} as the Zariski closure of the image of this map.
They derive a set of polynomial equations cutting out the rigid multiview variety and give a conjecture for a minimal generating set of its prime ideal.
The authors of \cite{joswig2016rigid} also describe defining polynomial equations for other variations of multiview varieties, e.g. when $\mathcal{X}$ in \eqref{eq:rigid-variation} is instead taken to be an arbitrary hypersurface in $\PP_\RR^3 \times \PP_\RR^3$, or when $\mathcal{X}$ is in $(\PP_\RR^3)^3$ and consists of triples of points with fixed pairwise Euclidean distances, or when $\mathcal{X}$ is in $(\PP_\RR^3)^4$ and consists of all $4$-tuples of points that lie on a common plane in $\PP_\RR^3$.

\textbf{Higher dimensions.}  Several works have studied higher-dimensional generalizations of the specialized joint camera map in \eqref{eq:jointImageMap}, i.e., maps from a projective space onto a product of projective spaces, all of arbitrary dimension, e.g. \cite{li2018images,ito2020projective,cid2021study}.
The study of higher-dimensional pinhole cameras is motivated by capturing \textit{dynamic scenes}.  We will sketch this application and some higher-dimensional extensions of multiview varieties in Section \ref{ssec:dynamic}.

\section{Understanding camera tuples} 
\label{sec:cameraTuples}
\subsection{Symmetry of joint camera map} The joint camera map $\Phi: \mathcal{X} \times \mathcal{C}_m \dashrightarrow \mathcal{Y}$ often carries a symmetry: depending on the specific camera models and 3D object types, there exists a group $G$ acting on the fibers of $\Phi$.  
  This symmetry arises by simultaneously acting on the cameras and the 3D objects via 3D change of coordinates which preserve the camera model and 3D objects. 
  This has the important implication that in 3D reconstruction (such as structure-from-motion) we cannot hope to precisely recover a world scene $(X, (C_1, \ldots, C_m))$ from its 2D image data $\Phi(X, (C_1, \ldots, C_m))$.  
  Instead the best we can reconstruct is the world scene up the action of $G$.  

For example, if $\mathcal{C}_m = (\PP\,\mathrm{Mat}^{3 \times 4}_3 )^m$ consists of tuples of projective pinhole cameras and $\mathcal{X} = (\PP_{\RR}^3)^n$ consists of tuples of points in 3D projective space, then the group $G$ is the projective general linear group $G = \operatorname{PGL}(4,\RR)$.  
Indeed, note 
\begin{align*}
   Ax = (Ag)(g^{-1}x)
\end{align*}
for each $A \in \mathrm{Mat}^{3 \times 4}_3$, $x \in \PP^3$ and $g \in \operatorname{PGL}(4,\RR)$.  In this case, we can only aim to recover the world scene up to a projective transformation.

As another example, suppose that that $\mathcal{C}_m$ consists of tuples of calibrated cameras $([R_1 | t_1], \ldots, [R_m | t_m])$
where $R_i \in \operatorname{SO}(3)$ and $t_i \in \RR^3$. 
Take $\mathcal{X} = (\PP^3)^{n}$ again to consist of tuples of 3D projective points.  Now the relevant group is 
\begin{align} \label{eq:similarity-group}
    G = \left\{ g \in \operatorname{GL}(4,\RR) : g = \begin{bmatrix} R & t \\ 0 & \lambda \end{bmatrix} \textup{ for some } R \in \operatorname{SO}(3), t \in \RR^3, \lambda \in \RR \setminus \{0\} \right\},
\end{align}
i.e., the scaled special Euclidean group of $\RR^3$ (also known as the similarity group).  This is because $G$ is the largest subgroup of 3D changes of coordinates which maps calibrated cameras to calibrated cameras.  
the world scene is ambiguous up to
a proper rigid motion $R$, and a central scale $\lambda$, $t$
So when we know the cameras are calibrated, the world scene is ambiguous up to a proper rigid motion ($R, t$) and a central scale ($\lambda$).  

A novel approach to detect all symmetries in a given 3D reconstruction problem has been recently developed by Duff, Korotynskiy, Pajdla, and Regan \cite{duff2022galois}.  Their approach is based on based on \textit{Galois / monodromy groups}, and uses numerical homotopy continuation computations (cf. Section~\ref{sec:homotopy-continuation} \nolinebreak below).

\subsection{Camera configurations} To cope with such ambiguities, it has been useful in computer vision to consider \textit{camera tuples modulo the relevant symmetries}.
Working out camera tuples modulo
symmetries can be viewed as one of the accomplishments of multiview geometry \cite{HZ}.
From a mathematical standpoint, Aholt and Oeding observed in \cite{aholt2014ideal} that formally we can regard the space of possible camera configurations as the GIT quotient:
\begin{align} \label{eq:GIT-quotient}
    \mathcal{C}_{m} / \!/ G.
\end{align}
Here GIT stands for geometric invariant theory, developed by Mumford \cite{mumford1994geometric}. 
The GIT quotient $\mathcal{C}_{m} /\!/ G$ is an abstract, non-embedded scheme or stack. 
While this GIT description can be made precise to the abstract algebraic geometer, it would not be useable in 3D reconstruction problems.
The perspective does tell us that multiview geometry found certain explicit birational models for these GIT quotients, which are useful for 3D reconstruction.

\subsection{Multifocal tensors} 
\label{ssec:cameraVarieties}

Camera configurations are represented by \textit{multifocal tensors} in multiview geometry.  
We present them for configurations of two, three and four cameras, and then give a higher-dimensional construction called \textit{Grassmann tensors} \cite{hartley2009reconstruction} motivated by dynamic scenes.

\subsubsection{\textup{\textbf{Two projective cameras: fundamental matrices}}} \label{ssec:fundamental}
Fix two projective pinhole cameras  $C_1, C_2 : \PP^3 \dashrightarrow \PP^2$ represented by matrices $A_1, A_2$.  
Consider the set of corresponding point pairs defined as the multiview variety:
\begin{equation*}
\mathcal{M}_{C} \, = \, \left\{ (x_1, x_2) \in \mathbb{P}^2 \times \mathbb{P}^2 : \exists X \in \mathbb{P}^3 \text{ s.t. } C_1(X) = x_1, C_2(X) = x_2 \right\}. 
\end{equation*}
This set is invariant under the action of $G = \operatorname{PGL}(4, \mathbb{R})$ on $C = (C_1, C_2)$, i.e. it depends only on the configuration $C$ modulo $G$.
On the other hand, $\mathcal{M}_{C}$ is generically a hypersurface of bidegree $(1,1)$ as explained in Example \ref{ex:fundamendalMatrix}: 
\begin{align*}
\mathcal{M}_C &= \bigg{\{} (x_1, x_2) : \exists X \in \mathbb{P}^3 \, \exists \lambda_1, \lambda_2 \in \mathbb{R} \text{ s.t. } \begin{bmatrix} A_1 & x_1 & 0 \\ A_2 & 0 & x_2 \end{bmatrix} \begin{small} \begin{bmatrix} X \\ -\lambda_1 \\ -\lambda_2 \end{bmatrix} \end{small} = 0 \bigg{\}} \\
& = \bigg{\{} (x_1, x_2 ) : \det \begin{small} \begin{bmatrix} A_1 & x_1 & 0 \\ A_2 & 0 & x_2 \end{bmatrix} \end{small} = 0 \bigg{\}}.
\end{align*}
This equation for $\mathcal{M}_C$ may be written as $x_2^{\top} F x_1 = 0$ where for $1 \leq i,j \leq 3$,
$$
F_{ij} = (-1)^{i+j} \det \begin{bmatrix} A_{1, \hat{j}} \\ A_{2, {\hat{i}}} \end{bmatrix},
$$
where $A_{1, \hat{j}}$ denotes $A_1$ with the $j$th row dropped and likewise for $A_{2, \hat{i}}$.
One regards $F$ as in $\mathbb{P}(\operatorname{Mat}^{3 \times 3})$, that is, defined only up to nonzero scale.   In multiview geometry, $F$ is called the \textit{fundamental matrix} for the pair of cameras $C$. 
It may be shown that $F$ generically determines $C$ modulo $G$ \cite{HZ}, thus $F$ furnishes an explicit model for the configuration.
The collection of all $F$ arising in this way is exactly the set of rank $2$ matrices, $\mathbb{P}(\operatorname{Mat}_2^{3 \times 3})$, with Zariski closure given by the determinantal variety $\{ F \in \mathbb{P}(\operatorname{Mat}^{3 \times 3}) : \det(F) = 0 \}$.
The construction is practically useful in 3D reconstruction because point pairs in $\mathcal{M}_C$ can be detected in real-life 2D images.

 \subsubsection{\textup{\textbf{Two calibrated cameras: Essential matrices}}} 
\label{ssec:essential}
 
 We can restrict the fundamental matrix construction to the case of two calibrated cameras.  
 The restriction yields a subset of the fundamental matrices, whose points are called \textit{essential matrices} and customarily denoted by $E$ in computer vision.
Demazure \cite{demazure1988deux} determined that the resulting subvariety of the determinantal variety is 5 dimensional, degree 10, with ideal minimally generated by the following ten cubic equations:
 \begin{equation}\label{eq:essential-eqns}
\det(E) = 0 \quad \quad \textup{and} \quad \quad E E^{\top} E - \frac{1}{2} \operatorname{trace}(EE^{\top})E = 0.
\end{equation}
 Another description is that essential matrices $E \in \mathbb{P}(\operatorname{Mat}^{3 \times 3})$ are those with rank $2$ whose top singular values are equal, i.e. $\sigma_1(E) = \sigma_2(E) > \sigma_3(E) = 0$.
 Yet another description is that, after taking the \textit{complex} Zariski closure, the variety of essential matrices can be realized as a complex hyperplane section of the determinantal variety of rank $\leq 2$ symmetric $4 \times 4$ matrices \cite{floystad2018chow}.
 
Here $E$ only depends on $C$ modulo the relevant group, i.e. the scaled special Euclidean group $G$ \eqref{eq:similarity-group}.  
In this case however, the association of a calibrated configuration to an essential matrix is generically two-to-one (cf. \textit{twisted pairs} in \cite{HZ}).  
Thus the GIT quotient \eqref{eq:GIT-quotient} gives a degree-2 cover of the variety of essential matrices. The advantage of essential matrices is that they capture corresponding image point pairs which can be detected in real-life images.

 \subsubsection{\textup{\textbf{Three cameras: Trifocal tensors}}}

Given three projective cameras $C_1, C_2, C_3 : \mathbb{P}^3 \dashrightarrow \mathbb{P}^2$ represented by matrices $A_1, A_2, A_3$, consider the set of corresponding point-line-line triples:
\begin{multline*}
\Psi_{C} = \big{\{} (x_1, \ell_2, \ell_3) \in \mathbb{P}^2 \times (\mathbb{P}^2)^*  \times (\mathbb{P}^2)^* : \\ \exists X \in \mathbb{P}^3 \text{ s.t. } C_1(X) = x_1, C_2(X) \in \ell_2, C_3(X) \in \ell_3 \big{\}},
\end{multline*}
where $(\mathbb{P}^2)^*$ is the dual projective plane consisting of lines in $\mathbb{P}^2$.  It then holds 
\begin{align*}
\Psi_C &= \bigg{\{} (x_1, \ell_2, \ell_3) : \exists X \in \mathbb{P}^3 \, \exists \lambda_1 \in \mathbb{R} \text{ s.t. } \begin{small} \begin{bmatrix} A_1 & x_1 \\[0.25em] \ell_2^{\top} A_2 & 0 \\[0.25em] \ell_3^{\top}A_3 & 0 \end{bmatrix} \begin{bmatrix} X \\ -\lambda_1 \end{bmatrix} \end{small} = 0  \bigg{\}} \\
& = \bigg{\{} (x_1, \ell_2, \ell_3) : \det \begin{small} \begin{bmatrix} A_1 & x_1 \\[0.25em] \ell_2^{\top} A_2 & 0 \\[0.25em] \ell_3^{\top}A_3 & 0 \end{bmatrix} \end{small} = 0 \bigg{\}},
\end{align*}
where we have identified $\ell_i$ with its normal vector of 3 homogeneous coordinates.
The  equation for $\Psi_C$ may be written using tensor contraction as $T(x_1, \ell_2, \ell_3) = 0$ where $T \in \mathbb{P}(\mathbb{R}^{3 \times 3 \times 3})$ is given by
  $$
 T_{ijk} \,\, = \,\, (-1)^{i+j+k} \, \det \begin{small} \begin{bmatrix} A_{1,i} \\ A_{2,j} \\ A_{3,\hat{k}} \end{bmatrix} \end{small}
 $$
 for $1 \leq i, j, k \leq 3$.
Here $A_{1,i}$ denotes the $i$th row of $A_1$, likewise for $A_{2,j}$, and $A_{3, \hat{k}}$ denotes $A_3$ with the $k$th row removed.
 In multiview geometry $T$ is known as the \textit{trifocal tensor} associated with the camera triple $C$. 
 The map from $(C_1,C_2,C_3)$ modulo $\mathrm{PGL}(4,\RR)$ to $T$ is generically one-to-one, so the Zariski closure of the set of trifocal tensors is a birational model for configurations of three (projective) cameras. 
 As proven in {\cite{aholt2014ideal}}, the trifocal tensors variety has dimension $18$, degree $297$ and an ideal minimally generated by $10$ polynomials in degree $3$, $81$ polynomials in degree $5$, and $1980$ polynomials in degree $6$.
 Trifocal tensors are relevant to 3D reconstruction, because of their relation to image correspondence data.
 They capture point-line-line correspondences by definition via tensor-vector contraction.  However it is also true that trifocal tensors capture other correspondence types (e.g., point-point-point) via other multilinear algebraic operations, see \cite[Part~III]{HZ} for details.

 If we restrict $(C_1,C_2,C_3)$ to  calibrated pinhole cameras, then we obtain only certain trifocal tensors and a subvariety called the \textit{calibrated trifocal variety} in \cite{kileel2017minimal}.  The calibrated trifocal variety is a birational model for the GIT quotient of triples of calibrated cameras modulo the scaled special Euclidean group $G$ \eqref{eq:similarity-group}.  
 Currently the ideal defining the variety is partially understood \cite{martyushev2017some}.

 \subsubsection{\textup{\textbf{Four cameras: Quadrifocal tensors}}}
A similar construction can be carried out for four projective pinhole cameras $C_1, C_2, C_3, C_4$ based on correspondences between quadruples of image lines:

\begin{multline*}
\big{\{}(\ell_1, \ell_2, \ell_3, \ell_4) \in (\mathbb{P}^2)^* \times (\mathbb{P}^2)^* \times (\mathbb{P}^2)^* \times (\mathbb{P}^2)^* : \\ \exists X \in \mathbb{P}^3 \text{ s.t. } C_1(X) \in \ell_1, C_2(X) \in \ell_2, C_3(X) \in \ell_3, C_4(X) \in \ell_4\big{\}}.
\end{multline*}
Similarly to the above, this is a hypersurface defined by $Q (\ell_1, \ell_2, \ell_3, \ell_4) = 0$ where $Q \in \mathbb{P}(\mathbb{R}^{3 \times 3 \times 3 \times 3})$ is given by
  $$
 Q_{ijkl} \,\, = \,\, (-1)^{i+j+k+l} \, \det  \begin{small} \begin{bmatrix}  A_{1,i} \\  A_{2,j}  \\  A_{3, k} \\ A_{4, l} \end{bmatrix} \end{small}
 $$
 The tensor $Q$ is known as the \textit{quadrifocal tensor} in multiple view geometry.  The ideal of the Zariski closure of quadrifocal tensors was studied up to degree $9$ in \cite{oeding2017quadrifocal}.
 Quadrifocal tensors for calibrated cameras have not been investigated.

\subsubsection{\textup{\textbf{Dynamic scenes}}} \label{ssec:dynamic}
The fundamental matrix and the trifocal and quadrifocal tensors (and their calibrated counterparts) are used for the reconstruction of static scenes. 
There are many approaches for the reconstruction of dynamic scenes, i.e., where the 3D points undergo different motions. 
Here we highlight two algebraic ideas.

Wolf and Shashua \cite{wolf2002projection} and Huang, Fossum and Ma \cite{huang2002generalized} argue that many motions of 3D points can be modeled by embedding them into a larger dimensional space such that the camera map becomes a projection of the form $\PP^N \dashrightarrow \PP^2$.
This motivated Hartley and Schaffalitzky \cite{hartley2009reconstruction} to construct the Grassmann tensor, which we explain in Section \ref{grassmannTensor} below.
For instance, a 3D point $X=(X_1,X_2,X_3,1)^\top$ that moves linearly in an affine chart of $\PP^3$ with constant velocity in the direction $V=(V_1,V_2,V_3,0)^\top$ can be embedded into $\PP^6$ as follows: $\bar X = (V_1:V_2:V_3:X_1:X_2:X_3:1)^\top$.
If $A_i$ denotes the $3 \times 4$ matrix representing a moving pinhole camera at time $i$, the image taken of the moving point $X$ at time $i$ is
$A_i (X+iV)$.
Writing $A'_i$ for the left-most $3 \times 3$ submatrix of $A_i$ and 
$\bar A_i = [i A'_i \mid A_i] \in \RR^{3 \times 7}$,
we can express the image as $A_i (X+iV) = \bar A_i \bar X$.
Hence, the camera map is a projection $\PP^6 \dashrightarrow \PP^2$.

Another algebraic approach to dynamic scenes was suggested by Vidal and coauthors \cite{vidal2006two,fan2006space}.
Suppose that a scene observed by a pinhole camera has $n$ independently moving objects and that we take two pictures of the scene at two different times.
For each of the $n$ objects, we have a fundamental matrix $F_i$ describing the relative pose of the object at the two instants of time.
Any pair of image points $(x_1,x_2)$ coming from a common 3D point must lie on one of the $n$ objects and thus satisfies the equation 
$(x_2^\top F_1 x_1) \cdots (x_2^\top F_n x_1)=0$.
As this equation is homogeneous of degree $n$ in both $x_1$ and $x_2$, 
we can apply the Veronese embedding $\nu_n: \PP^2 \to  \PP^N, N = \binom{n+2}{2}-1,$
to obtain the following bilinear equation in $\nu(x_1)$ and $\nu(x_2)$:
$\nu(x_2)^\top F \nu(x_1) = 0$,
where $F \in \RR^{N \times N}$ is called the \emph{multibody fundamental matrix}.
This construction has been generalized, e.g. to \emph{multibody trifocal tensors} \cite{hartley2004multibody} and other motion models \cite{vidal2006unified}.
Finally, we note that 3D motion segmentation from $m$ views, i.e., the problem to determine the number $n$ of moving objects and which given $m$-tuples of image points belong to which object, can also be solved with \emph{generalized principal component analysis (GPCA)}. GPCA is a method considered by Vidal, Ma, and Sastry to fit an unknown number of linear subspaces of unknown dimensions to given data points in a real vector space \cite{vidal2005generalized}.
Other algebraic methods for GPCA include \cite{tsakiris2017filtrated}.

 \subsubsection{\textup{\textbf{Grassmann tensors}}} \label{grassmannTensor}
 Here we review the Grassmann tensor construction of Hartley and Schaffalitzky \cite{hartley2009reconstruction} that unifies multifocal tensors with higher-dimensional generalizations. 
Consider $\mathbb{P}^N$, noting that $N > 3$ can model dynamic scenes as in Section~\ref{ssec:dynamic}.
Let  $C_i : \mathbb{P}^N \dashrightarrow \mathbb{P}^{n_i}$, $X \mapsto A_i X$ be $m$ surjective linear projections.  Assume $m \leq N+1 \leq \sum_{i=1}^m n_i$. 
Then we can fix integers $c_1, \ldots, c_m$ satisfying $1 \leq c_i \leq n_i$ and $c_1 + \ldots + c_m = N+1$. 
{The special cases $n_i = 2$ and $n_i = 1$ appear in pinhole cameras and 1D radial cameras \cite{larsson2020calibration,thirthala2012radial}, respectively.}

Write $\Gr(n-c, \mathbb{P}^{n})$ for the Grassmannian of linear subspaces of $\mathbb{P}^n$ of codimension $c$.
The main object is the set of ``corresponding subspace tuples":
\begin{multline} \label{eq:setup-Grassmann}
\{(L_1, \ldots, L_m) \in \Gr(n_1-c_1, \mathbb{P}^{n_1}) \times \ldots \times \Gr(n_m - c_m, \mathbb{P}^{n_m}) : \\ \exists X \in \mathbb{P}^N \text{ s.t. } C_1(X) \in L_1, \ldots, C_m(X) \in L_m \}.
\end{multline}
The authors of \cite{hartley2009reconstruction} proved that \eqref{eq:setup-Grassmann} is an algebraic hypersurface; moreover its defining equation is multilinear in the Pl\"ucker coordinates for the Grassmannians. 
Thus, the coefficients of defining equation form an $m$-way tensor
$$
\mathbf{T} = \mathbf{T}_{(C_1, \ldots, C_m), (c_1, \ldots,  c_m)}
$$
of size $\binom{n_1+1}{c_1} \times \ldots \times \binom{n_m+1}{c_m}$.  One calls $\mathbf{T}$ the  \emph{Grassmann tensor} for the projections $(C_1, \ldots, C_m)$ with respect to the codimensions $(c_1, \ldots, c_m)$  \cite{hartley2009reconstruction} . 
The authors proved that $\mathbf{T}$ depends only on $(C_1, \ldots, C_m)$ modulo $\operatorname{PGL}(N+1, \mathbb{R})$.  Further, if at least $n_i$ strictly exceeds $1$ then $(C_1, \ldots, C_m)$ modulo $\operatorname{PGL}(N+1)$ is generically uniquely determined by $\mathbf{T}$.  
Fundamental matrices, trifocal tensors and quadrifocal tensors correspond to $N=3$ and $(c_1, c_2) = (2,2)$, $(c_1, c_2, c_3) = (2,1,1)$ and $(c_1, c_2, c_3, c_4) = (1,1,1,1)$ respectively.

The construction of $\mathbf{T}$ can be derived here. 
For $L_i  \in \Gr(n_i-c_i, \mathbb{P}^{n_i})$  we choose a basis of the dual space $L_{i}^{*} \subseteq (\mathbb{P}^{n_i})^\ast$ comprised of the hyperplanes containing $L_i$; let the basis be
\begin{equation*}
\mathcal{V}_i = \{v_{i1}, \ldots, v_{ic_i}\} \subseteq (\mathbb{P}^{n_i})^{\ast}.
\end{equation*}
The dual linear map $C_i^{*} : (\mathbb{P}^{n_i})^{\ast} \rightarrow (\mathbb{P}^N)^{\ast}$ pulls $\mathcal{V}_i$ back to $(\mathbb{P}^N)^{\ast}$ to give
\begin{equation*}
 A_i^{\top} \mathcal{V}_i = \{A_i^{\top}v_{i1}, \ldots, A_i^{\top}v_{ic_i}\} \subseteq (\mathbb{P}^N)^{\ast}.
\end{equation*}
The requirement on $X \in \mathbb{P}^N$ in \eqref{eq:setup-Grassmann} is  $\langle A_iX, v_{ij} \rangle =0$ for each $i,j$, or equivalently $\langle X, A_i^{\top} v_{ij} \rangle =0$.
Hence we require $A_1^{\top}\mathcal{V}_1 \cup \ldots \cup A_m^{\top} \mathcal{V}_{m}$ to not span $(\mathbb{P}^N)^{\ast}$, i.e. 
\begin{equation} \label{eq:my-mat}
    \begin{bmatrix}
    A_1^{\top} v_{11} & \ldots & A_1^{\top} v_{1c_1} & \ldots & \ldots & A_{m}^{\top} v_{m1} & \ldots & A_{m}^{\top} v_{m c_m} \\
    \end{bmatrix} \in \mathbb{R}^{(N+1) \times (N+1)}.
\end{equation}
is rank-deficient.  The matrix \eqref{eq:my-mat} factors as
\begin{equation} \label{eq:rect_prod}
\underbrace{\begin{bmatrix}
A_1^{\top} & A_2^{\top} & \ldots & A_m^{\top} \\
\end{bmatrix}}_{= A^{\top}}
\underbrace{\begin{bmatrix}
\mathcal{V}_1 & 0 & \ldots & 0 \\
0 & \mathcal{V}_2 & \ldots & 0 \\
\vdots & \vdots & \ddots  & \vdots \\
0 & 0 & \ldots & \mathcal{V}_m
\end{bmatrix}}_{= \mathcal{V}}.
\end{equation}
We require the determinant of \eqref{eq:rect_prod} to vanish, which by the Cauchy-Binet formula is
\begin{align} \label{eq:Cauchy-Binet}
    \sum_{\substack{\mathcal{I} \subseteq [(n_1+1) + \ldots + (n_m + 1)]\\ \# \mathcal{I} = N+1}} \det \! \left( A^{\top}  [:, \mathcal{I}] \right) \, \det \! \left(\mathcal{V}[\mathcal{I},:] \right). 
\end{align}
Here $A^{\top} [:,\mathcal{I}]$ denotes the $(N+1) \times (N+1)$ submatrix  consisting of the columns indexed by $\mathcal{I}$, etc.  
Recalling that the rows in $\mathcal{V}$ containing $\mathcal{V}_i \in \mathbb{R}^{(n_i+1) \times c_i}$ have rank $c_i$ and that $c_1 + \ldots + c_m = N+1$, we can restrict the sum in \eqref{eq:Cauchy-Binet} as follows:
\begin{equation*}
    \sum_{\substack{\mathcal{I}_1 \subseteq [n_1 + 1] \\ \# \mathcal{I}_1 = c_1}} \!\! \ldots \!\! \sum_{\substack{\mathcal{I}_m \subseteq [n_m + 1] \\ \# \mathcal{I}_m = c_m}} \!\! \det \left( A_1^{\top}[:,\mathcal{I}_1] \, \ldots \, A_m^{\top}[:, \mathcal{I}_m] \right) \, 
    \prod_{i=1}^m \det\left(\mathcal{V}_i[\mathcal{I}_i,:] \right).
\end{equation*}

 At this point, note $\det\left(\mathcal{V}_i[\mathcal{I}_i,:] \right)$ for $\mathcal{I}_i \subseteq [n_i+1]$ are the \textup{dual Pl\"ucker coordinates} of $L_i \in \Gr(n_i-c_i, \mathbb{P}^{n_i})$. 
These depend on the choice of basis $\mathcal{V}_i$ only through a global nonzero scale, and (after sign flips) equal the primal Pl\"ucker coordinates of $L_i$ up to scale.
That is, letting $\mathcal{U}_i = \{u_{i1}, \ldots, u_{i,n_i-c_i+1}\} \subseteq \mathbb{P}^{n_i}$ be a basis for $L_i$ one has
\begin{multline} \label{eq:prop}
\Big{(}\det(\mathcal{V}_i[\mathcal{I}_i,:]) \;\;\mid\;\; \mathcal{I}_i \subseteq [n_i+1], \# \mathcal{I}_i = c_i \Big{)} \\ \propto \,\,\, \Big{(}(-1)^{c_i + \sum \mathcal{I}_i }\det(\mathcal{U}_i[\mathcal{I}_i^c,:]) \;\;\mid\;\; \mathcal{I}_i \subseteq [n_i+1], \# \mathcal{I}_i = c_i \Big{)},
\end{multline}
where $\mathcal{I}_i^c = [n_i+1] \setminus \mathcal{I}_i$ and $\sum \mathcal{I}_i = \sum_{\ell \in \mathcal{I}_i} \ell$.

In view of this discussion, we define the entries of the Grassmann tensor as 
\begin{equation} \label{eq:grassmann-tensor}
    \mathbf{T}_{\mathcal{I}_1, \ldots, \mathcal{I}_m} \,\, = \,\, (-1)^{N+1 + \sum \mathcal{I}_1 + \ldots + \sum \mathcal{I}_m} \,\, \det \! \left( A_1^{\top}[:,\mathcal{I}_1] \, \ldots \, A_m^{\top}[:, \mathcal{I}_m] \right),
\end{equation}
and conclude the following result:

\begin{theorem}[\cite{hartley2009reconstruction}]
Given surjective linear projections $C_i : \mathbb{P}^N \dashrightarrow \mathbb{P}^{n_i}$ for $i = 1, \ldots, m$ and
  $c_i \in \mathbb{Z}$ such that $1 \leq c_i \leq n_i$ and $c_1 + \ldots + c_m = N+1$,  
define the set of corresponding subspaces by \eqref{eq:setup-Grassmann}. 
Up to Zariski closure, this is an algebraic hypersurface cut out in primal Pl\"ucker coordinates by  $\mathbf{T}$ in \eqref{eq:grassmann-tensor}.
\end{theorem}

 \subsubsection{\textup{\textbf{Beyond four cameras}}}  \label{sec:beyond4}
When $N=3$ and there are $m>4$ pinhole cameras, a tensor-based model for camera configurations is no longer possible since there is no suitable tuple of codimensions $(c_1, \ldots, c_m)$ for the Grassmann tensor construction.

The issue is that various embedded subvarieties that one may associate with the cameras are not hypersurfaces, and thus not specified by a single equation.
 However models for camera configurations still exist inside appropriate \textit{Hilbert schemes}.  
 Roughly speaking, a Hilbert scheme is the set of ideals of a fixed polynomial ring with a prescribed number of equations in each degree.
 In particular, one can associate $(C_1, \ldots, C_m)$ modulo $G$ to the ideal of its multiview variety $\mathcal{M}_C$, viewed as a point in a suitable Hilbert scheme \cite{aholt2013hilbert}. 
 Aholt, Thomas and Sturmfels \cite{aholt2013hilbert} proved that this association is generically one-to-one if $m>2$.  
 The Hilbert scheme approach was recently extended by Lieblich and Van Meter \cite{lieblich2020two} via functorial methods.  Their work has the benefit of incorporating calibrated cameras as well.
 
A rather different line of work encodes the configuration of $m$ cameras through their \emph{viewing graph}.
 The graph has $m$ vertices, one for each  camera.
 There is an edge between two vertices if the relative position of the two corresponding cameras is known.
 The edge is then labeled with the fundamental matrix of the camera pair.
 In most practical scenarios with many cameras, the viewing graph is incomplete.
 The structure of the graph then determines if the given fundamental matrices are generically sufficient to recover the global camera configuration
\cite{trager2018solvability,arrigoni2021viewing}.

\section{3D reconstruction} \label{sec:3Dalgorithms}
In this section, we discuss the role that multiview geometry plays in 3D reconstruction algorithms and how zero-dimensional polynomial systems arise. The next section will focus on methods that are used to solve these polynomials with sufficient speed.

\subsection{Structure-from-motion pipeline} \label{ssec:pipeline}

\begin{figure}[htb]
    \centering
    \includegraphics[width=\textwidth]{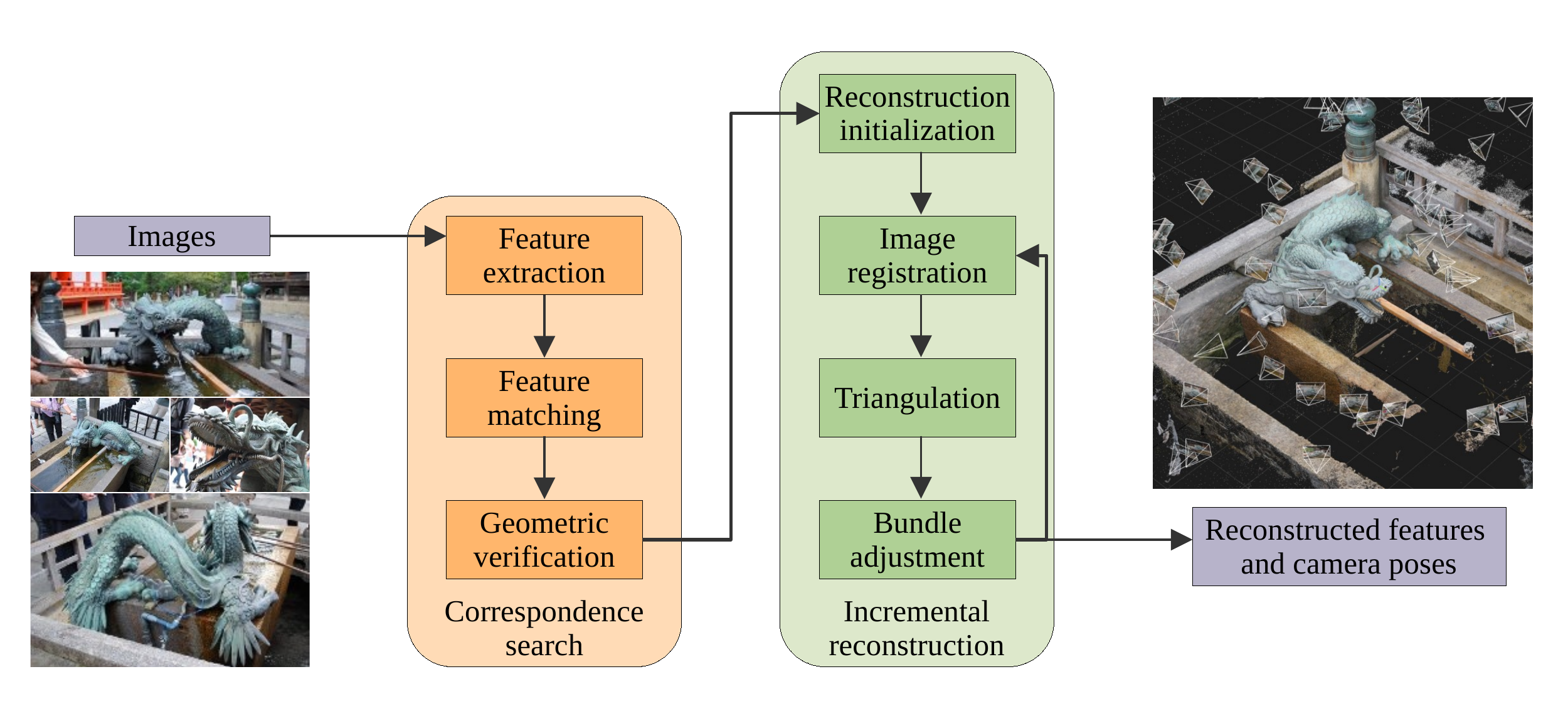}
    \caption{Incremental structure-from-motion pipeline.}
    \label{fig:pipelineDiagram}
\end{figure}

We have introduced the structure-from-motion problem  already in Section~\ref{sec:computervisionproblems}.
Here, we give more detailed algorithmic steps to solve this 3D reconstruction problem.  
Figure~\ref{fig:pipelineDiagram} overviews the usual \emph{incremental structure-from-motion pipeline}; see \cite{bianco2018evaluating} for a full description.

\textbf{Feature extraction.}
The first step in the pipeline is to detect distinguished points or lines (or possibly other information such as conics) in each given image. 
A standard algorithm is the Scale-Invariant Feature Transform (SIFT) due to Lowe \cite{lowe2004distinctive}.

\textbf{Feature matching.}
Next, the detected features are matched across different images; cf. Figure~\ref{fig:matching}. 
In particular, feature matching identifies which images show a common part of the scene.

\textbf{Geometric verification.}
Typically, not all feature matches from the feature matching step will be correct; many will be outliers.
For two (or more) images with matched features, the geometric verification step tries to compute a geometric transformation that correctly maps sufficiently many features between the images.
If such a transformation can be found, the pair (or tuple) of images and their common features that could be mapped are considered to be geometrically verified.
For instance, the fundamental (see Section \ref{ssec:fundamental}) or essential matrix (see Section \ref{ssec:essential}) encode the relative pose of two pinhole cameras and thus provide the desired transformation between their two images.
The output of the geometric verification step is the viewing graph whose nodes are images and whose (hyper-)edges connect geometrically verified image pairs (or tuples).
The (hyper-)edges are labeled with the computed geometric transformations.
See Section~\ref{sec:beyond4} for more on the viewing graph.

Since geometric verification forgets the original images and only works with the feature matches from feature matching (see Figure \ref{fig:algebraStep}), its core is an algebro-geometric problem that amounts to solving systems of polynomial equations. 
However, due to the presence of outliers, this step must use robust estimation techniques, usually RANSAC (Random Sample Consensus).
More details are in Section \ref{ssec:ransac}.

\textbf{Reconstruction initialization.}
This step picks a starting pair (or tuple) of geometrically verified images (typically, in a dense area of the viewing graph), marks them as \emph{registered}, and reconstructs the 3D coordinates of their matched features (as shown in Figure \ref{fig:algebraStep}).
Note that the relative pose of the cameras was already precomputed in the previous step.

\textbf{Image registration.}
This is the first step of the incremental reconstruction phase.
It registers one new image that is connected with previously registered images in the viewing graph
by (re-)calculating its camera pose relative to the reconstructed camera poses of the registered images.
This can be either done from matched 2D features (as in the geometric verification step) or from 2D-3D correspondences if some of the already reconstructed 3D features are detected 2D features on the new image. 
Both ways solve a system of polynomial equations, but the latter problem -- known as \emph{Perspective-$n$-Point (PnP)} if all detected features are points -- can be typically solved faster. 

\textbf{Triangulation.}
This step reconstructs the 3D coordinates of all features that the newly registered image has in common with any of the previously registered images (except, of course, for those features whose 3D coordinates were reconstructed before); see Section \ref{ssec:edd}.

\textbf{Bundle adjustment.}
This step refines the 3D features and camera poses reconstructed so far to minimize the accumulation of errors in the incremental reconstruction phase; see \cite{triggs1999bundle} for more details.
  Afterward, a new image will be registered and the incremental reconstruction phase repeats until (almost) all cameras and sufficiently many 3D features are reconstructed to approximate the 3D scene.
  Bundle adjustment is a very important step in the pipeline, but a fuller account falls outside our scope.  See \cite{ozyecsil2017survey} for another survey article that emphasizes numerical optimization aspects in structure-from-motion more ours does.

\subsection{Geometric verification} \label{ssec:ransac}

The algebro-geometric problems within the structure-from-motion pipeline aim to reconstruct the relative camera positions locally.  That is, we seek the relative pose for only a few (e.g., two) of the given images at a time, based on feature matches between the images.
A salient issue is that some of the purported feature matches will be \textit{completely erroneous}.
Indeed feature extraction/matching methods applied to real images invariably produce a nontrivial fraction of incorrect matches; and the problem becomes severe for scenes with repeated textures (e.g., a brick wall), among others.  Therefore reconstruction methods must be robust to outliers.  A general paradigm for estimation in the presence of outliers was introduced by Fischler and Bolles \cite{fischler1981random}, who were motivated by the computer vision setting.  Their \textit{Random Sample Consensus (RANSAC)} and its subsequent variants remain popular today for computing local reconstructions.

To illustrate RANSAC and show how it motivates solving zero-dimensional polynomial systems in vision, we discuss essential matrix estimation. In this example, data consists of several purported point matches 
$$(x_1, y_1), \ldots, (x_n, y_n) \in \mathbb{P}_{\mathbb{R}}^2 \times \mathbb{P}_{\mathbb{R}}^2.$$   
Let $\mathcal{E} \subset \mathbb{P}(\operatorname{Mat}^{3 \times 3})$ be the variety of essential matrices defined by \eqref{eq:essential-eqns}.
We seek $E \in \mathcal{E}$ which best fits the veridical point pairs amongst the data.
By construction, if $(x_i, y_i)$ is a correct (and noiseless) image point pair and $E$ is the ground-truth essential matrix, then 
$$
 y_i^{\top} E x_i = 0.
$$
Notice that this constraint on $E$ is linear and of codimension $1$.   
Thus given $\dim(\mathcal{E}) = 5$ such constraints, we expect up to $\operatorname{deg}(\mathcal{E}) = 10$  real solutions for $E$.
However it is unknown which of the purported points pairs are veridical.
RANSAC achieves robust estimation by randomly sampling a subset $\mathcal{I} \subseteq [n] = \{1, \ldots, n\}$ of five point pairs and fitting to them exactly by solving the polynomial system:
\begin{equation} \label{eq:solveE}
 \det(E) = 0,  \quad \quad  E E^{\top} E - \frac{1}{2} \operatorname{trace}(EE^{\top})E = 0, \quad \quad    y_{i}^{\top} E  x_i = 0 \,\,\, (i \in \mathcal{I}).
\end{equation}
Letting $E_1, \ldots, E_s$ be the real solutions, RANSAC checks if there exists a consensus amongst the other $n-5$ point pairs as to whether any of the solutions agree with most of the other point pairs. 
RANSAC re-samples $\mathcal{I}$ until a consensus is found. 
We note that this involves a choice of error metric to measure the agreement between $E_k$ and $(x_i, y_i)$ for $i \in \mathcal{I}$.  
A simple choice is $(y_i^{\top} E_k x_i)^2$ often called algebraic error, but often the geometric reprojection error or its first-order approximation called  Sampson error are preferred \cite{HZ}.  
We also remark that popular variants of RANSAC, like LO-RANSAC \cite{chum2003locally}, polish the winning $E_k$ by using all point pairs in $[n]$ that agree with $E_k$ to form a non-linear least squares cost which is minimized starting from $E_k$.

\subsection{Minimal problems}
\label{sec:minimalProblems}

The polynomial system \eqref{eq:solveE} is called a \textup{minimal problem} in the vision literature, reflecting the fact that \eqref{eq:solveE} uses the minimal amount of data such that the solution is uniquely determined up to finitely many solutions.
More generally, minimal problems are 3D reconstruction problems such that random input instances have a finite positive number of solutions.
Formally, a structure-from-motion problem is \emph{minimal} if the associated joint camera map $\Phi: \mathcal X \times \mathcal C_m \dashrightarrow \mathcal Y$ over the complex numbers has generically finite non-empty fibers, modulo the group $G$ that acts on the fibers.
In other words, for a generic (complex) $y \in \mathcal Y$, the fiber $\Phi^{-1}(y)$ is non-empty and consists of finitely many  $G$-orbits.

By choosing appropriate coordinates for camera configurations modulo the group $G$ (for example, by using the multifocal tensors in Section~\ref{ssec:cameraVarieties}), a minimal problem can be cast as a zero-dimensional parameterized polynomial system.
Thus one expects that a minimal problem can be solved efficiently, as long as its generic number of complex solutions is sufficiently small.
This generic number of complex solutions -- known as the \emph{algebraic degree} of the minimal problem -- is a measure of the minimal problem's intrinsic difficulty.

\iffalse
A first indicator if a sufficiently fast minimal solver for a (camera-)minimal problem can be developed is the \emph{algebraic degree} of the problem, i.e., its number of complex solutions $N_{\CC}$ for generic image tuples.
This degree measures the intrinsic difficulty of the problem.
For instance, the algebraic degrees of the (camera-)minimal problems in Examples \ref{ex:5ptProblem} and \ref{ex:cameraMinimal} are $20$ and $272$, respectively.
For a given real image tuple $y \in \mathcal{Y}$, the number of real solutions $N_{\RR, y}$ depends on the data $y$.
%
\textcolor{blue}{
However, some integers $n$ appear with positive probability as the number of real solutions;
more specifically, we call $n \in \mathbb{Z}_{\geq 0}$ a \textit{typical real root count} if the subset $\{ y \in \mathcal{Y} \mid N_{\RR, y} = n \} $ has positive Lebesgue measure in the data space $\mathcal{Y}$.
Often, there are multiple typical real root counts, and they are hard to determine.  However, each must be less than or equal to $N_{\CC}$, and congruent to $N_{\CC}$ modulo \nolinebreak $2$.
}
Since minimal solvers typically first find all complex solutions and afterwards only keep the real ones, the larger the algebraic degree, the more difficult it is to develop a sufficiently fast minimal solver.
\fi

\begin{example} \label{ex:5ptProblem}
Reconstructing five points in general position observed by two calibrated pinhole cameras (see Figure \ref{fig:algebraStep}) is a minimal problem.
This is because, generically, the system \eqref{eq:solveE} has $10$ complex solutions (i.e., essential matrices), so $20$ complex pairs of cameras (modulo the scaled special Euclidean group), and for each camera pair there are five unique 3D points mapping to the given image points. 
Note that the number of \textit{real} solutions to \eqref{eq:solveE} depends on the image data.  The real count is generically amongst $\{0,2,4,6,8,10\}$, but varies with the specific five points.
In particular, understanding when there exists at least a real solution is subtle \cite{agarwal2017existence}.
An even more difficult problem is to determine when there is at least one chiral reconstruction, i.e. a real solution where the 3D points lie \emph{in front} of the cameras \cite{pryhuber2020existence}.
If we are given six points in two calibrated views, the system \eqref{eq:solveE} no longer corresponds to a minimal problem, as generically there is no complex solution.  We need the image data to lie on the \textit{Chow hypersurface} in order for  \eqref{eq:solveE} to be soluble; see \cite{floystad2018chow} for an explicit Pfaffian formula defining the Chow hypersurface obtained using Ulrich sheaves, or \cite{bik2020jordan} for a simpler construction based on Jordan algebras.
\end{example}

The literature on minimal problems is vast.
Many minimal problems have been described, solvers have been developed for them, and new minimal problems are constantly appearing; see e.g. \cite{johansson2002structure,oskarsson2004minimal,kukelova2007minimal,nister2007minimal,stewenius2008minimal,ramalingam2008minimal,elqursh2011line,kuang2013pose,kukelova2015efficient,ventura2015efficient,saurer2015minimal,lee2016minimal,camposeco2016minimal,barath2017minimal,barath2018five,miraldo2018minimal,zhao2019minimal,kukelova2017clever,fabbri2019trplp}. 

There are also many variations of the above definition.
For instance, when solving 3D reconstruction by RANSAC as described in Section~\ref{ssec:ransac}, one often aims to recover the camera parameters only. 
Hence, often a structure-from-motion problem is referred to as minimal if for a generic image tuple $y \in \mathcal Y$, the projection of the fiber $\Phi^{-1}(y) \subset X \times \mathcal C_m$ onto the space $\mathcal C_m$ of cameras is non-empty and finite modulo $G$.
In the following, we call such problems \emph{camera-minimal} and reserve the term \emph{minimal} for when we aim to recover the 3D structure in $\mathcal X$ as well (as in the definition before Example \ref{ex:5ptProblem}).
Note that every minimal problem is camera-minimal but not vice versa. 

\begin{example} \label{ex:cameraMinimal}
The following problem is camera-minimal (shown in \cite{kileel2017minimal} where it is called $3\mathrm{PPP}+1\mathrm{PPL}$): Reconstruct four points $X_1,X_2,X_3,X_4 \in \PP^3$ and a line $L \subset \PP^3$ incident to $X_4$ that are observed by three calibrated pinhole cameras such that the first two cameras see only the four points and the third camera only sees $X_1,X_2,X_3$ and the line $L$; see Figure \ref{fig:joeminimal}.
In fact, it typically has $272$ complex solutions in terms of camera parameters.
However, this problem is \emph{not} minimal since the line $L$ cannot be recovered uniquely when the three cameras are known:
There is a one-dimensional family of lines that yield the same image line in the third view.
\begin{figure}[htb]
    \centering
    \includegraphics[width=0.8\textwidth]{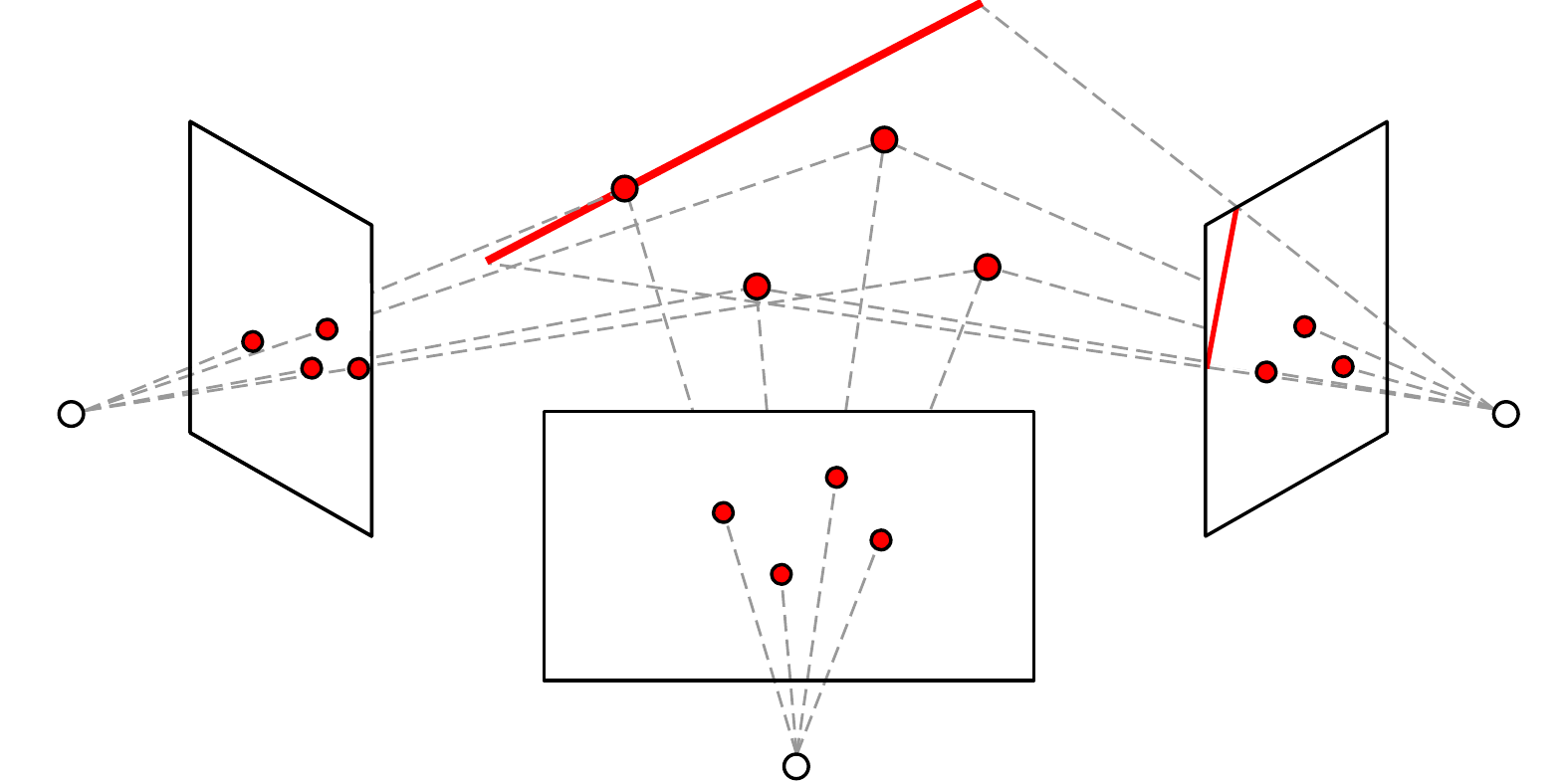}
    \caption{A camera-minimal problem with algebraic degree $272$.}
    \label{fig:joeminimal}
\end{figure}
\end{example}

In standard structure-from-motion reconstructions, (camera-)minimal problems need to be solved many times to ensure RANSAC finds a solution that is not corrupted by outliers. 
Thus, minimal solvers need to be fast.  
Furthermore, they should be accurate, since numerical errors could jeopardize RANSAC's success as well. 
We also would like to handle various nearly-degenerate geometries, such as near coplanarities amongst the five point pairs in \eqref{eq:solveE}, as such situations do arise in natural images.  
In Section~\ref{ssec:eliminationTemplates}, we review the symbolic and numerical minimal solvers that are used to meet  these needs.

The large zoo of minimal problems in the literature is motivated by the quest to find reconstruction problems that 
1) admit a feasible minimal solver as described above, 
2) their necessary input data can be detected by state-of-the-art feature extraction methods, and
3) cover the wide variety of reconstruction applications in computer vision.
The latter includes variations of the camera model (see Section \ref{sec:cameramodels}), of the parameters to be reconstructed, or of the input data (also beyond the structure-from-motion problem we focus on here, e.g. camera pose estimation from known 3D \emph{and} 2D data).
We list some of the most classical minimal problems in Table~\ref{tab:minimal}.

\begin{table}[htb]
    \centering
\scriptsize
 \begin{tabular}{|| c | c |  c | c||} 
 \hline
  {{\textbf{to be reconstructed}}} & {{\textbf{minimal data}}} & {{\textbf{ degree}}} & {{\textbf{references}}} \\  
 \hline\hline
essential matrix & 5 point pairs & 10  & \cite{nister2004efficient,demazure1988deux}  \\
 \hline
\centered{fundamental matrix} & \centered{7 point pairs} & \centered{3} &  \centered{\cite{hartley1994projective} \\ already known in 19th \\   century \cite{hesse1863cubische,sturm1869problem}} \\
 \hline
 \centered{relative pose of 2 calibrated \\ cameras with  unknown  \\ common  focal length }
 & \centered{6 point pairs} & \centered{15}  & \centered{\cite{stewenius2008minimal}}    \\
 \hline 
 \centered{absolute pose of 1  \\
calibrated camera  \\ (P3P, image registration)}
  & \centered{3 world-image  point pairs }&  \centered{4} 
  & \centered{\cite{haralick1994review}   \\ already known in 19th \\ century \cite{grunert}}\\
 \hline
planar homography & 4 point pairs & 1 &  \\ \hline
calibrated planar homography & 4 point pairs & 12 & \cite{duff2019plmp,faugeras1988motion} \\
 \hline
trifocal tensor & 9 line triples & 36 & \cite{larsson2017efficient}  \\ \hline
\centered{calibrated trifocal tensor} 
& \centered{3 point triples \\ +1 line triple}
& \centered{216} & \centered{\cite{kileel2017minimal}} \\
\hline
\centered{relative pose of 2 projective  \\
cameras with unkown \\
radial lens distortion}
& \centered{8 point pairs} & \centered{16} & \centered{\cite{kukelova2007minimal}}  \\
\hline 
world point under noise
 & known cameras with: && \\
(triangulation, & $\bullet$ 1 point pair \hspace{0.9mm} & $\bullet$ 6\hspace{1.4mm} & \cite{hartley1997triangulation} \\
reprojection error)
&$\bullet$ 1 point triple & $\bullet$ 47 &
\cite{stewenius2005hard} \\ \hline
\end{tabular}
    \caption{Some representative minimal problems.  Several of these minimal problems are classical and integrated in vision pipelines.}
    \label{tab:minimal}
\end{table}

Efforts have been undertaken to provide complete catalogs of some classes of minimal problems.
Kileel \cite{kileel2017minimal} describes all $66$ camera-minimal problems (incl. their algebraic degrees) for three calibrated pinhole cameras that can be formulated using linear constraints on the trifocal tensor.
Duff, Kohn, Leykin, and Pajdla \cite{duff2019plmp} characterize all $30$ minimal problems (including their algebraic degrees) of reconstructing points, lines and their incidences for an arbitrary number of calibrated pinhole cameras under the assumption that all points and lines are observed by every camera.
In a follow-up work allowing partial visibility \cite{duff2020pl}, they classify all camera-minimal and minimal problems of three calibrated pinhole cameras observing point-line arrangements such that each line is incident to at most one point.
In particular, there are $759$ (equivalence classes of) camera-minimal problems with algebraic degree less than $300$.

\section{Polynomial solvers for minimal problems} \label{ssec:eliminationTemplates}

In this section, we discuss symbolic and numerical methods to construct solvers for minimal problems in computer vision, e.g. the minimal problems in Table~\ref{tab:minimal}. 

Algebraically, we model a minimal problem as a \textit{zero-dimensional parameterized polynomial system:}
\begin{equation} \label{eq:param-system}
\begin{cases} 
    f_1(x_1, \ldots, x_n; \theta_1, \ldots, \theta_t) = 0 \\
    \hspace{2cm} \vdots \\
    f_k(x_1, \ldots, x_n; \theta_1, \ldots, \theta_t) = 0.
\end{cases}
\end{equation}
Here $x=(x_1, \ldots, x_n)$ are variables, and $\theta = (\theta_1, \ldots, \theta_t)$ are parameters (e.g., coming from image point correspondences).
The equations $f_i$ are polynomials in both $x$ and $\theta$, typically with integer coefficients (e.g., the epipolar constraints \eqref{eq:solveE}).
Write 
   $\mathcal{I}$
for the ideal generated by the equations
\begin{equation*}
\mathcal{I} = \{ p_1 f_1 + \ldots + p_k f_k \,\mid\, p_i \in \mathbb{R}[x_1, \ldots, x_n; \theta_1, \ldots, \theta_t]\}.
\end{equation*}
Let $|_{\theta}$ denotes specialization to concrete parameter values $\theta \in \mathbb{R}^t$.
Minimality of the problem translates to the following two assumptions:
\begin{itemize}
    \item[\textbf{A1)}] \textup{For generic $\theta$, the ideal $\mathcal{I}|_{\theta}$ is zero-dimensional (i.e., has finitely many roots)};
    \item[\textbf{A2)}] \textup{For generic $\theta$, the ideal $\mathcal{I}|_{\theta}$ is radical (i.e., has no repeated roots)}
\end{itemize}
The conditions imply the existence of a positive integer $N_{\CC} = N_{\CC}(\mathcal{I})$ such that for generic parameters $\theta$, the number of $\mathbb{C}$-solutions to \eqref{eq:param-system} is $N_{\CC}$.  Furthermore, each root has multiplicity $1$. 
Note that $N_{\mathbb{C}}$ is the algebraic degree of the minimal problem  discussed in Section~\ref{sec:minimalProblems}.

By definition, a \textit{solver} is meant to receive values for the parameter $\theta \in \mathbb{R}^t$, and then output the set of corresponding real solutions $V_{\mathbb{R}}(\mathcal{I}|_{\theta}) \subset \mathbb{R}^n$ to \eqref{eq:param-system}.
A salient issue to bear in mind is the \textit{{online vs. offline}} distinction in computer vision:  
When given a specific parameter value $\theta$, for many vision tasks we must run the solver in \textup{real-time}.  
Therefore, online computations must be fast.
On the other hand, when developing the solver, we are free to dedicate far more time and computational resources offline.

\subsection{Symbolic method: Elimination templates} \label{sec:elim-temp}

Currently a symbolic method is the most popular one to build minimal solvers in vision.  
At the heart of the method is a classical fact that zero-dimensional polynomial systems can be converted to eigenvector problems.\footnote{A different symbolic approach, based instead on \textit{resultants}, is being developed by Bhayani, Heikkil\"a and Kukelova  \cite{bhayani2023sparse,bhayani2020sparse}.}

\begin{lemma}\label{lem:eigenval}
Let $\theta \in \mathbb{R}^t$, and the corresponding solutions to \eqref{eq:param-system} be $V_{\mathcal{\mathbb{C}}}(\mathcal{I}|_{\theta}) = \{v_1, \ldots, v_{N_{\mathbb{C}}}\} \subset \mathbb{C}^n$.
Let $P \in \mathbb{R}[x_1, \ldots, x_n]$ be a polynomial separating the solutions, that is, $P(v_i) \neq P(v_j)$ for all $i \neq j$.
Let $\mathcal{B} = \{b_1, \ldots, b_{N_{\mathbb{C}}} \} \subseteq \mathbb{R}[x_1, \ldots, x_n]$ be such that its reduction modulo $\mathcal{I}|_{\theta}$ is a vector space basis $\overline{\mathcal{B}} = \{\overline{b}_1, \ldots, \overline{b}_{N_{\mathbb{C}}} \}$ of the quotient ring $\mathbb{R}[x_1, \ldots, x_n] / \mathcal{I}|_{\theta}$.
Multiplication by $P$ is a linear endomorphism on $\mathbb{R}[x_1, \ldots, x_n] / \mathcal{I}|_{\theta}$, represented with respect to $\overline{\mathcal{B}}$ by a matrix
\begin{equation} \label{eq:mult-p-def}
    [\operatorname{mult}(P)|_{\theta}]_{\overline{\mathcal{B}}} \, \in \, \mathbb{R}^{N_{\CC} \times N_{\CC}},  \quad  \quad  \overline{P} \overline{b}_j  := \sum_{i=1}^{N_{\CC}} \left([\operatorname{mult}(P)|_{\theta}]_{\overline{\mathcal{B}}}\right)_{ij} \overline{b}_i.
\end{equation}
Then $[\operatorname{mult}(P)|_{\theta}]_{\overline{\mathcal{B}}}$ is diagonalizable with eigenvalue/left eigenvectors pairs given by $(P(v_i), \mathcal{B}(v_i)) \in \mathbb{C} \times \mathbb{C}^{N_{\mathbb{C}}}$ for $i = 1, \ldots, N_{\mathbb{C}}$.
\end{lemma}
In \cite{cox2021stickelberger}, Cox gives an interesting history of Lemma~\ref{lem:eigenval}, which is often (wrongly) attributed to Stickelberger.
We call $P$ the \textit{action polynomial} and $[\operatorname{mult}(P)|_{\theta}]_{\overline{\mathcal{B}}}$ the \textit{action matrix}.
Notice that if the action matrix is available, then Lemma~\ref{lem:eigenval} shows how the solutions $v_i$ to \eqref{eq:param-system} may be computed: from the eigenvectors $w_i$ of the action matrix we require $\mathcal{B}(v_i)$ to be proportional to $w_i$.
To make the recovery of $v_i$ easy, typically $\mathcal{B}$ is chosen to consist of monomials.  

Thus the main challenge is to obtain the action matrix online. 
In principle, it requires computing identities of the form 
\begin{equation} \label{eq:q-identity}
    Pb_j = \sum_{i=1}^{N_{\CC}} \left([\operatorname{mult}(P)|_{\theta}]_{\overline{\mathcal{B}}}\right)_{ij} b_i \,\, + \,\, \sum_{i=1}^k q_{ij} f_i|_{\theta} \,\, \in \,\, \mathbb{R}[x_1, \ldots, x_n],
\end{equation}
for $j = 1, \ldots, N_{\mathbb{C}}$ where $q_{ij} \in \mathbb{R}[x_1, \ldots, x_n]$ and $\left([\operatorname{mult}(P)|_{\theta}]_{\overline{\mathcal{B}}}\right)_{ij} \in \mathbb{R}$.
Gr\"obner basis computations are too slow to be performed online, and not viable in floating point arithmetic anyways.  Instead, the main tool is \textit{elimination templates}.

Fix choices of $\mathcal{B} \subset \mathbb{R}[x_1, \ldots, x_n]$ and $P \in \mathbb{R}[x_1, \ldots, x_n]$ such that: 
\begin{enumerate}
\item[\textbf{A3)}] \textup{For generic $\theta$,  $\overline{\mathcal{B}} \subset \mathbb{R}[x_1, \ldots, x_n]/\mathcal{I}|_{\theta}$ is a vector space basis;}  
    \item[\textbf{A4)}] \textup{For generic $\theta$, $P$ separates $V_{\mathbb{C}}(\mathcal{I}|_{
    \theta})$.} 
\end{enumerate}
In the present setting, elimination templates are defined as follows.

\begin{definition}\label{def:elim-templ}
A set $\mathcal{E} = \{e_1, \ldots, e_T\} \subset \mathcal{I} \subseteq \mathbb{R}[x_1, \ldots, x_n; \theta_1, \ldots, \theta_t]$ is an \emph{elimination template} if  \textbf{A5)} For generic $\theta$, we have $P \cdot \mathcal{B} \subset \operatorname{span}(\mathcal{B} \cup \mathcal{E}|_{\theta})$.
\end{definition} 
\noindent Often $\mathcal{E}$ is chosen to be $x$-monomials multiplied by $f_1, \ldots, f_k$.  
By Definition~\ref{def:elim-templ},
\begin{equation} \label{eq:with-star}
    Pb_j = \sum_{i=1}^{N_{\CC}} \left([\operatorname{mult}(P)|_{\theta}]_{\overline{\mathcal{B}}}\right)_{ij} b_i \,\, + \,\, \sum_{\ell=1}^T \star \, e_{\ell}|_{\theta} \,\, \in \,\, \mathbb{R}[x_1, \ldots, x_n]
\end{equation}
for $j =1, \ldots, N_{\mathbb{C}}$ and scalars $\star$. 
Computing the action matrix thus becomes a linear system. Specifically, fix $\mathcal{C} = \{c_1, \ldots, c_{M} \} \subset \mathbb{R}[x_1, \ldots, x_n]$ such that 
\begin{enumerate}
\item[\textbf{A6)}] {For generic $\theta$,  $\mathcal{E}|_{\theta} \subset \operatorname{span}(\mathcal{B} \cup \mathcal{C})$.}  
\end{enumerate}
Often $\mathcal{C}$ is all $x$-monomials appearing in any $e_{\ell} \in \mathcal{E}$ with nonzero coefficient in $\mathbb{R}[\theta_1, \ldots, \theta_t]$,  not already in $\mathcal{B}$.
Form a matrix called the \textit{elimination template matrix}, denoted  $[\operatorname{elim}(P)|_{\theta}]_{\mathcal{B},\mathcal{E},\mathcal{C}} 
\in \mathbb{R}^{(T+{N_{\CC}}) \times ({N_{\CC}} + M)}$, with rows indexed by $\mathcal{E}|_{\theta} \, \cup \, P \cdot \mathcal{B}$ and columns indexed by $\mathcal{B} \cup \mathcal{C}$:

\begin{equation} \label{eq:template-matrix}
\begin{pNiceArray}{ccc|ccc}[first-row,first-col,margin]
    & & \begin{matrix} b_1 & \ldots & b_{N_{\CC}} \end{matrix} &  &  & \begin{matrix} c_1 & \ldots & c_{M} \end{matrix} & \\
\begin{matrix} e_1|_{\theta} \\ \vdots \\ e_T|_{\theta} \end{matrix} & & \begin{matrix} \\ \scalebox{3.5}{$\star$} \\ \end{matrix} & & & \begin{matrix} \\ \scalebox{3.5}{$\star$} \\ \end{matrix}    &  \\
& & & & & & \\\hline
\begin{matrix} Pb_1 \\ \vdots \\ Pb_{N_{\CC}} \end{matrix} & & \begin{matrix} \\ \scalebox{3.5}{$\star$} \\ \end{matrix} & & & \begin{matrix} \\ \scalebox{3.5}{$\star$} \\ \end{matrix} 
\end{pNiceArray},
\end{equation}
where each entry is the coefficient of the column in the row. 
Thanks to \eqref{eq:with-star}, elementary row operations can reduce the elimination template matrix to the form

\vspace{-0.2em}

\begin{equation} \label{eq:reduce-matrix}
\begin{pNiceArray}{ccc|ccc}[first-row,first-col,margin]
    & & \begin{matrix} b_1 & \ldots & b_{N_{\CC}} \end{matrix} &  &  & \begin{matrix} c_1 & \ldots & c_{M} \end{matrix} & \\
\begin{matrix} e_1|_{\theta} \\ \vdots \\ e_T|_{\theta} \end{matrix} & & \begin{matrix} \\ \scalebox{3.5}{$\star$} \\ \end{matrix} & & & \begin{matrix} \\ \scalebox{3.5}{$\star$} \\ \end{matrix}    &  \\
& & & & & & \\\hline
\begin{matrix} Pb_1 \\ \vdots \\ Pb_{N_{\CC}} \end{matrix} & & \begin{matrix} \\ \scalebox{1.3}{$[\operatorname{mult}(P)|_{\theta}]_{\overline{\mathcal{B}}}$} \\ \end{matrix} & & & \begin{matrix} \\ \scalebox{2}{$0$} \\ \end{matrix} 
\end{pNiceArray}.
\end{equation}
The bottom-left block of \eqref{eq:reduce-matrix} is the desired action matrix.

The eigenvectors and elimination templates method is summarized in Algorithm~\ref{alg:solver}.
In the online part, the most costly step is computing the pseudo-inverse  of the top-right block of the elimination template matrix. %, taking $\mathcal{O}(TM(T+M))$ flops.  
In the offline phase, we seek $\mathcal{B}, P, \mathcal{E}, \mathcal{C}$ such that the size of the elimination template matrix is as small as possible.

\begin{algorithm} 
  \caption{Offline/online solver for the parameterized system \eqref{eq:param-system}} \label{alg:solver}
  \begin{algorithmic}[1]
    \OfflineInputs{$\langle f_1, \ldots, f_k\rangle  \subseteq \mathbb{R}[x_1, \ldots, x_n; \theta_1, \ldots, \theta_t]$ satisfying \textbf{A1}-\textbf{A2}. \\[10pt]} 
    \OfflinePre{Find $\mathcal{B} \subseteq \mathbb{R}[x_1, \ldots, x_n]$ satisfying \textbf{A3}, \\ $P \in \mathbb{R}[x_1, \ldots, x_n]$ satisfying \textbf{A4}, \\ $\mathcal{E} \subseteq \mathbb{R}[x_1, \ldots, x_n; \theta_1, \ldots, \theta_t]$ satisfying \textbf{A5}, \\ $\mathcal{C} \subseteq \mathbb{R}[x_1, \ldots, x_n]$ satisfying \textbf{A6}. \\[7pt]
    Compute the bottom ${N_{\CC}}$ rows of \eqref{eq:template-matrix}, and \\
    denote these by $\begin{bmatrix}
    E_{21} & E_{22}
    \end{bmatrix}  \in \mathbb{R}^{{N_{\CC}} \times ({N_{\CC}}+M)}$. \\[10pt]} 
    \OnlineInputs{Generic $\theta \in \mathbb{R}^t$. \\[10pt]} 
    \OnlineComp{
    Evaluate the top $T$ rows of \eqref{eq:template-matrix} at $\theta$, and \\
    denote these by $\begin{bmatrix}
    E_{11} & E_{12}
    \end{bmatrix} \in \mathbb{R}^{T \times ({N_{\CC}}+M)}.$ \\
    Thus $[\operatorname{elim}(P)|_{\theta}]_{\mathcal{B}, \mathcal{E}, \mathcal{C}} = \begin{bmatrix}
E_{11} & E_{12} \\
E_{21} & E_{22}
\end{bmatrix}$.
    \\[7pt]
    Compute $E_{21} - E_{22}E_{12}^{\dagger}E_{11} \in \mathbb{R}^{{N_{\CC}} \times {N_{\CC}}}$.\\
    This equals $[\operatorname{mult}(P)|_{\theta}]_{\overline{\mathcal{B}}}$. \\[7pt]
    Compute all real eigenvalue/eigenvector pairs $(\lambda_1,w_1), \ldots, (\lambda_r,w_r) \in \mathbb{R} \times \mathbb{R}^{N_{\CC}}$ of $[\operatorname{mult}(P)|_{\theta}]_{\overline{\mathcal{B}}}$. \\[7pt]
    Determine $v_1, \ldots, v_r \in \mathbb{C}^n$ such that \\ $P(v_i) = \lambda_i$ and $\mathcal{B}(v_i) \propto w_i$.
    \\[10pt]} 
    \OnlineOut{$V_{\mathbb{R}}(\mathcal{I}|_{\theta}) = \{v_1, \ldots, v_r\} \cap \mathbb{R}^n$.}
  \end{algorithmic}
\end{algorithm}

For the offline procedure, full details are beyond our scope but we mention some of the highlights.  A breakthrough was obtained by Kukelova, Bujnak and Pajdla in 2008, who showed how find elimination templates automatically, as opposed to by hand \cite{kukelova2008automatic}.  Behind their work, a principle is that ``shape" of symbolic calculations should be the same at generic parameters $\theta$, and also preserved if we reduce \eqref{eq:param-system} modulo a prime.  
This enables traditional Gr\"obner basis computations in the offline phase.  
In \cite{larsson2017efficient} Larsson, \AA{}str\"om and Oskarsson observed that $q_{ij}$ in \eqref{eq:q-identity} are uniquely defined only up to the first \textup{syzygy module} of  $(f_1|_{\theta}, \ldots, f_k|_{\theta})$. 
Syzygy-reduction \cite[Proposition~1]{larsson2017efficient} was shown to result in choices of $\mathcal{E}$ with substantially smaller elimination template matrices for many minimal problems.  
In \cite{larsson2018beyond}, Larsson, Oskarsson, \AA{}str\"om, Wallis, Kukelova and Pajdla leveraged the fact that the basis $\mathcal{B}$ need not be standard monomials with respect to some Gr\"obner basis for $\mathcal{I}|_{\theta}$,  only condition \textbf{A3} is required (see \cite{mourrain1999new}).
The authors chose $\mathcal{B}$ by a random monomial-sampling strategy, and showed the strategy can sometimes produce smaller elimination template matrices than achievable anywhere on the Gr\"obner fan. 
Martyushev, Vrablikova and Pajdla further improved the state-of-the-art with a greedy search strategy to reduce the size of $\mathcal{C}$ in \cite{martyushev2022optimizing}. 
Integrating many of these ideas, Li and Larsson released code for building solvers in \cite{li2020gaps}.

\subsection{Numerical method: Homotopy continuation}\label{sec:homotopy-continuation}
An alternative method for building minimal problem solver is homotopy continuation. 
 Homotopy continuation \cite{sommese2005numerical} is a path-following, ODE-based numerical method to solve parametrized zero-dimensional polynomial systems of a significantly different nature than the eigenvectors and elimination templates method described above.
It uses only floating-point arithmetic operations, and exploits the offline vs. online paradigm in a different way.

In the offline phase, homotopy continuation solves one generic instance of the minimal problem (i.e., the polynomial system $f(x; \theta_*) = 0$ for a fixed choice of parameters).  
The offline solve is done effectively from scratch.  
Specifically, homotopy continuation chooses a \textit{start system} 
$g(x) = 0$
of polynomials with the same variables and degree structure as $f(x; \theta_*) = 0$ that is trivial to solve (e.g., the start solutions are roots of unity).
Then homotopy continuation interpolates from the start system to the minimal problem instance via linear interpolation: 
$$h(x; t) = (1-t) f(x; \theta_*) + tg(x),$$ 
as $t$ goes from $1$ to $0$.  During the interpolation, we numerically track how the solutions to $g(x) = 0$ evolve \cite{davidenko1953new,davidenko1953approximate}. 
In general, many solution paths diverge.  However the paths with finite limits give all solutions to $f(x; \theta_*) = 0$ \cite[Theorem~7.1.1]{sommese2005numerical}. 
The offline procedure is costly: its complexity is driven by the number of solutions to $g(x) = 0$, often exponential in the number of variables.

In the online phase, homotopy continuation also interpolates between polynomial systems. 
But now it tracks solutions of the precomputed minimal problem instance $f(x; \theta_*) = 0$ to the solutions of the desired instance of the minimal problem (i.e., $f(x; \theta) = 0$ for a different choice of parameters). 
The complexity of the online phase is driven by the number of solutions to $f(x; \theta_0) =0$, rather than the number of solutions to $g(x) =0$.
Consequently the online solve is much faster than the offline phase.
Indeed, in the online solve we track only $N_{\mathbb{C}}$ many paths, which is the intrinsic algebraic degree of the minimal problem.
General-purpose software for homotopy continuation includes Bertini \cite{bates2013numerically}, HomotopyContinuation.jl \cite{breiding2018homotopycontinuation},  NumericalImplicitization.m2 \cite{chen2019numerical}, among others.

\begin{figure}[htb]
    \centering
    \includegraphics[width=0.6\textwidth]{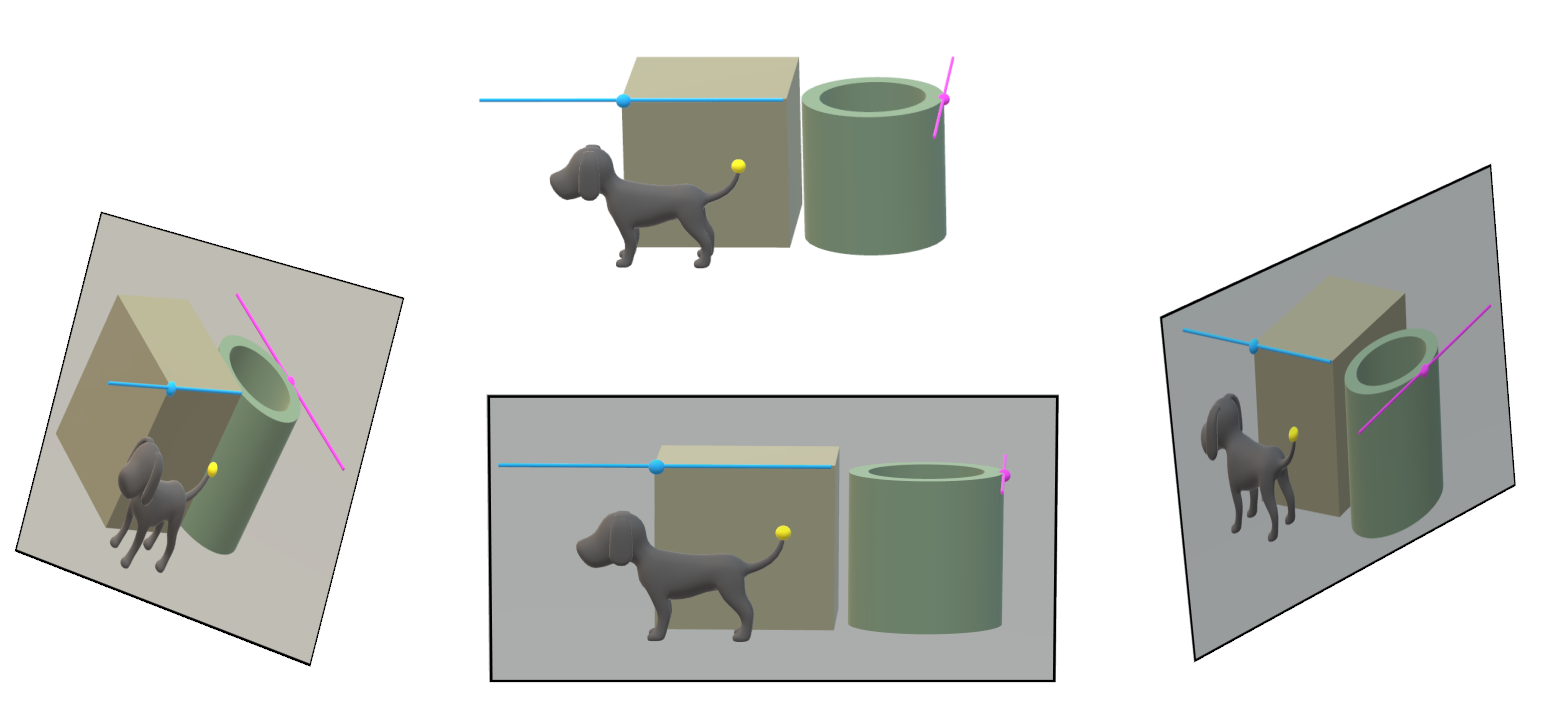}
    \caption{A minimal problem with algebraic degree $312$: three calibrated pinhole cameras observe three points, two of which have an incident line.  Homotopy continuation was used to build a solver for this minimal problem taking $1$ second or less, whereas symbolic methods haven't yielded a practical solver \cite{fabbri2019trplp}.}
    \label{fig:minimal312}
\end{figure}

In computer vision, homotopy continuation has been used to determine the algebraic degrees of (camera-)minimal problems \cite{kileel2017minimal,duff2019plmp,duff2020pl}.
Homotopy continuation was also previously suggested as a solver for such problems \cite{horn1991relative,kriegman1992geometric}, however for many years it was not seen to be competitive with the elimination templates approach.
Recently, homotopy
continuation yielded the fastest solver for the minimal structure-from-motion problem shown in Figure~\ref{fig:minimal312}.
For this problem no  elimination template matrix of reasonable size could be found by the symbolic method of Section~\ref{sec:elim-temp}. 
However, with homotopy continuation a numerically stable solver was built that runs on the order of a second or less \cite{fabbri2019trplp}.   The minimal problem in Figure \ref{fig:minimal312} has algebraic degree $N_{\mathbb{C}} = 312$, so homotopy continuation tracks $312$ paths per online solve. 
The solver was validated on real-data scenarios where the state-of-the-art pipeline
COLMAP \cite{schonberger2016structure} failed to find enough point-based features,
which also highlights the importance of investigating (camera-)minimal problems that involve not only points but also lines.

More recently, homotopy continuation has been combined with machine learning to build competitive solvers for minimal problems of large algebraic degree: 
Hruby, Duff, Leykin, and Pajdla \cite{hruby2022learning} achieve fast computational speeds by first using machine learning to pick a single start solution and then tracking that solution to a single solution of the desired target system.
The start solution does not come from a fixed start system, but rather determines the start system used for homotopy continuation, meaning that the start system is also learned. 
Another important advance has been to combine GPU computing with homotopy continuation.  This was recently carried out by Chien, Fan, Abdelfattah,  Tsigaridas, Tomov and Kimia  \cite{chien2022gpu} on various benchmark minimal problems.

\section{Degeneracies and discriminants} \label{sec:discriminants}

This section is devoted to degenerate configurations of 3D scenes and cameras.
Section~\ref{ssec:critical} focuses on configurations that do not admit unique recovery from image data.
Degenerate views of 3D curves and surfaces are discussed in Section~\ref{ssec:aspect}.
Finally, Section~\ref{ssec:hurwitz} outlines the connection between degenerate instances of minimal problems and condition numbers that measure how sensitive numerical reconstruction algorithms are to errors in the input.

\subsection{Critical loci} 
\label{ssec:critical}

A \emph{critical configuration} in a structure-from-motion problem is a pair $(X,C) \in \mathcal{X}\times\mathcal{C}_m$ consisting of an ordered subset  $X$ of $\PP_\RR^3$ and an $m$-tuple of cameras $C$ such that the resulting $m$ images $\Phi(X,C)$ do not have a unique reconstruction modulo the group $G$ that acts on the fibers of the joint camera map $\Phi$.
In other words, a configuration of 3D points $X \in \mathcal{X}$ and cameras $C \in \mathcal{C}_m$ is critical if there is another configuration $(X',C') \in \mathcal{X}\times\mathcal{C}_m$ that is not in the same $G$-orbit as $(X,C)$ and produces the same images, i.e., $\Phi(X,C)=\Phi(X',C')$.

\begin{figure}[htb]
    \centering
    \includegraphics[width=0.8\textwidth]{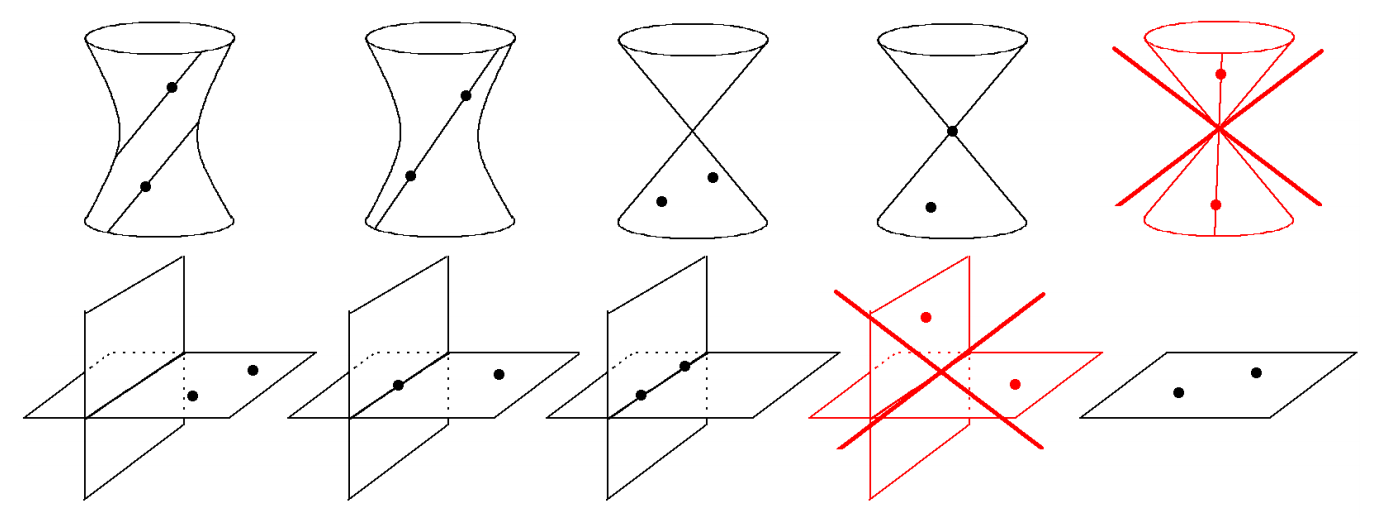}
    \caption{The eight critical quadrics (in black) for 3D points and two projective pinhole cameras.
    The marked points are the camera centers.
    The exceptions in Theorem \ref{thm:critical} are red. (Courtesy of Martin Bråtelund.)}
    \label{fig:critical}
\end{figure}

The study of critical configurations for two pinhole cameras goes back to the German-Austrian literature on photogrammetry, where critical loci where known as \emph{gefährliche Örter}:
Krames \cite{krames1941ermittlung} showed that an ordered set of 3D points and two pinhole cameras form a critical configuration only if the points and the two camera centers lie on a real ruled quadric surface; see also Maybank \cite{maybank2012theory} for a modern reference.
Recall that a quadric surface in $\PP_\RR^3$ is said to be \emph{ruled} if it contains a line. 
There are four projectively non-equivalent ruled quadrics.
They can be distinguished based on the rank of a symmetric matrix $S \in \RR^{4 \times 4}$ giving a defining equation $X^\top S X = 0$ of the quadric: 
smooth ruled quadric ($\mathrm{rank}\,S= 4$), cone ($\mathrm{rank}\,S= 3$), two planes ($\mathrm{rank}\,S= 2$), and double plane ($\mathrm{rank}\,S= 1$); see Figure~\ref{fig:critical}.

A first characterization of  critical configurations of 3D points and an \emph{arbitrary} number of projective pinhole cameras modulo the group $\mathrm{PGL}(4, \RR)$ was provided by Hartley and Kahl \cite{hartley2007critical}.
In particular, for $m=2$, they derive:

\begin{theorem} \label{thm:critical}
An ordered set of 3D points and two pinhole cameras form a critical configuration modulo $\mathrm{PGL}(4, \RR)$ if and only if the points and the two camera centers lie on a real ruled quadric surface $Q$, with the following two exceptions:
1) $Q$ is a cone and both camera centers lie on the same line on $Q$ but not on its vertex, and 
2) $Q$ is two planes and the camera centers do not lie on the same plane.
Hence, up to projective equivalence, there are eight critical quadrics with two marked camera centers, shown in Figure~\ref{fig:critical}.
\end{theorem}

The classification in \cite{hartley2007critical} has some mistakes for three or more cameras. Bråtelund developed new algebraic techniques \cite{braatelund2021critical2} and corrected the classification for three cameras \cite{braatelund2021critical3}. 

Kahl and Hartley \cite{kahl2002critical} also studied the critical configurations of 3D points and calibrated pinhole cameras, modulo the scaled special Euclidean group of $\RR^3$.
They provide a full characterization for two cameras and show the existence of critical configurations for an arbitrary number of points and cameras.
It is still an open problem to derive the complete catalog of all critical configurations consisting of 3D points and more than two calibrated pinhole cameras.

\begin{remark}
Another natural notion of critical locus that is not as prominent in the vision literature as the critical configurations described above is the ramification locus of the joint camera map. 
For instance, for a minimal problem of algebraic degree larger one, almost all configurations $(X,C) \in \mathcal{X}\times\mathcal{C}_m$ are critical (over $\mathbb{C}$).
However, the ramification locus of the joint camera map of such a minimal problem carries crucial meaning: configurations close to the ramification locus often correspond to ill-conditioned  problem instances \cite{burgisser2013condition,demmel1987condition,demmel1987geometry}, which has negative algorithmic consequences such as numerical instability; 
see also Section \ref{ssec:hurwitz}.
\end{remark}

Critical configurations are also understood in a variety of other settings.
Critical loci for 3D lines and three pinhole cameras are for instance investigated by Buchanan \cite{buchanan1992criticalConf,buchanan1992critical}, Maybank \cite{maybank1995critical}, Navab and Faugeras \cite{navab1997critical}, and Zhao and Chung \cite{zhao2008critical}.
They are congruences in $\Gr(1,\PP^3)$ that can be parametrized by Bordiga surfaces in $\PP^4$. 
It would be interesting to generalize those results to more cameras, which is currently open.

Buchanan \cite{buchanan1988twisted} also described the critical configurations for a single pinhole camera, i.e., for the problem of \emph{camera calibration} where both 3D and image points are known and only the camera parameters are to be recovered.
Here the critical 3D points and the camera center lie either on a (possibly reducible but connected) twisted cubic or on the union of a line and a plane such that the camera center lies on the line.

Åström and Kahl \cite{aastrom2003ambiguous} classified the critical loci for 1D structure-from-motion problems where 2D points are projected via pinhole cameras of the form $\PP^2 \dashrightarrow \PP^1$.
Critical configurations of 2D points and camera centers lie on cubic curves.

Bertolini, Turrini, and coauthors generalized several of the results above to higher dimensions, i.e., where the classical pinhole cameras $\PP^3 \dashrightarrow \PP^2$ are replaced by projections $\PP^N \dashrightarrow \PP^n, N > n$; see e.g. \cite{bertolini2007critical,bertolini2015critical,bertolini2020critical,bertolini2020criticalB,bertolini2021smooth}.
As explained at the end of Section \ref{ssec:cameraVarieties}, the study of general projections is motivated by the reconstruction of dynamic scenes.
Several of their works also perform simulated experiments to test numerical instability phenomena for problem instances near critical loci; see in addition \cite{bertolini2007instability,bertolini2007instabilityB,bertolini2019critical}.

Despite the literature on critical configurations outlined above, critical loci of vision problems are still unexplored in many settings, e.g. for 3D point-line arrangements or for other camera models besides pinhole cameras.

\subsection{Aspect graphs and visual events} \label{ssec:aspect}
\begin{figure}[htb]
    \centering
    \includegraphics[width = 0.4\textwidth]{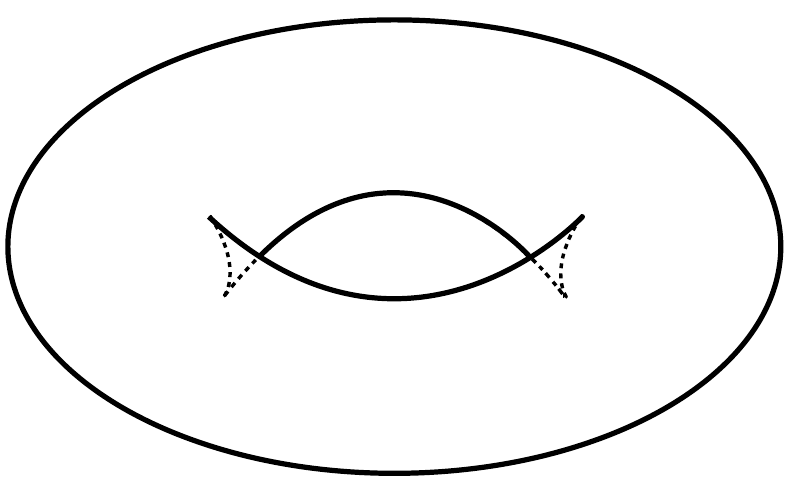}
    \caption{Image contour of a torus.}
    \label{fig:contour}
\end{figure}

In this section, we focus on taking pictures of an algebraic curve or surface $X$ in $\PP^3$ with a pinhole camera $C: \PP^3 \dashrightarrow \PP^2$. 
Clearly, if $X$ is a curve, its image is typically a plane curve and we write $Y_C(X)$ for its Zariski closure in $\PP^2$.
If $X$ is a surface, we denote by $Y_C(X)$ its \emph{image contour}, also known as \emph{silhouette} or \emph{outline curve}, in the image plane $\PP^2$.
This is the curve that bounds the two-dimensional region that is the image of the surface $X$ taken by the camera $C$.
In other words, it is the natural sketch one might use to depict the surface; see Figure \ref{fig:contour} for an image contour of a torus.
Formally, $Y_C(X)$ is the branch locus of the map $C$ restricted to $X$, i.e., it is the projection of the critical points where viewing lines through the camera center are tangent to the surface $X$.

The reconstruction of an algebraic surface from its image contour has been studied intensively, e.g., by Zariski \cite{zariski1929problem}, Segre \cite{segre},  Biggiogero \cite{biggiogero1947caratterizzazione,biggiogero1947sulla}, Chisini \cite{chisini1947sulla} who conjectured that a surface can be reconstructed from its image contour, Kulikov \cite{kulikov1999chisini,kulikov2008chisini} who solved Chisini's conjecture, and Forsyth \cite{forsyth1996recognizing} who solved the special case of a smooth surface and a general camera. 
Algorithmic solutions are for instance provided in \cite{kriegman1990recognizing,d1992courve,gallet2019reconstruction,gallet2021reconstruction}.

The \emph{aspects} of a curve or surface $X$ are the stable views when the camera moves, i.e., where the topology and singularities of the plane curve $Y_C(X)$ do not change under perturbations of the camera center. 
There is a finite number of aspects and they are the nodes of the \emph{aspect graph}.
Its edges are the \emph{visual events} that transition between different aspects. 
For instance, Figure \ref{fig:contour} shows an aspect of a torus, and a visual event occurs when we rotate the torus until the real singularities disappear.
The aspect graph was first introduced by Koenderink and van Doorn \cite{koenderink1979internal} under the name of \emph{visual potential}.
A detailed discussion of image contours and their visual events can be found in Koenderink's book  \cite{koenderink1990solid}.

Although there has not been much use of aspect graphs in real-life applications, they have been an active research topic in the computer vision community; see e.g. Bowyer's and Dyer's survey \cite{bowyer1990aspect} or Chapter~13 in Forsyth's and Ponce's book \cite{forsyth2011computer}.
In particular, several algorithms for the computation of aspect graphs of algebraic surfaces were proposed, using both symbolic and numerical methods \cite{ponce1990computing,petitjean1992computing,rieger1992global,rieger1993computing,pae2001computing}.

\begin{figure}[htb]
    \centering
    \includegraphics[width = 0.8\textwidth]{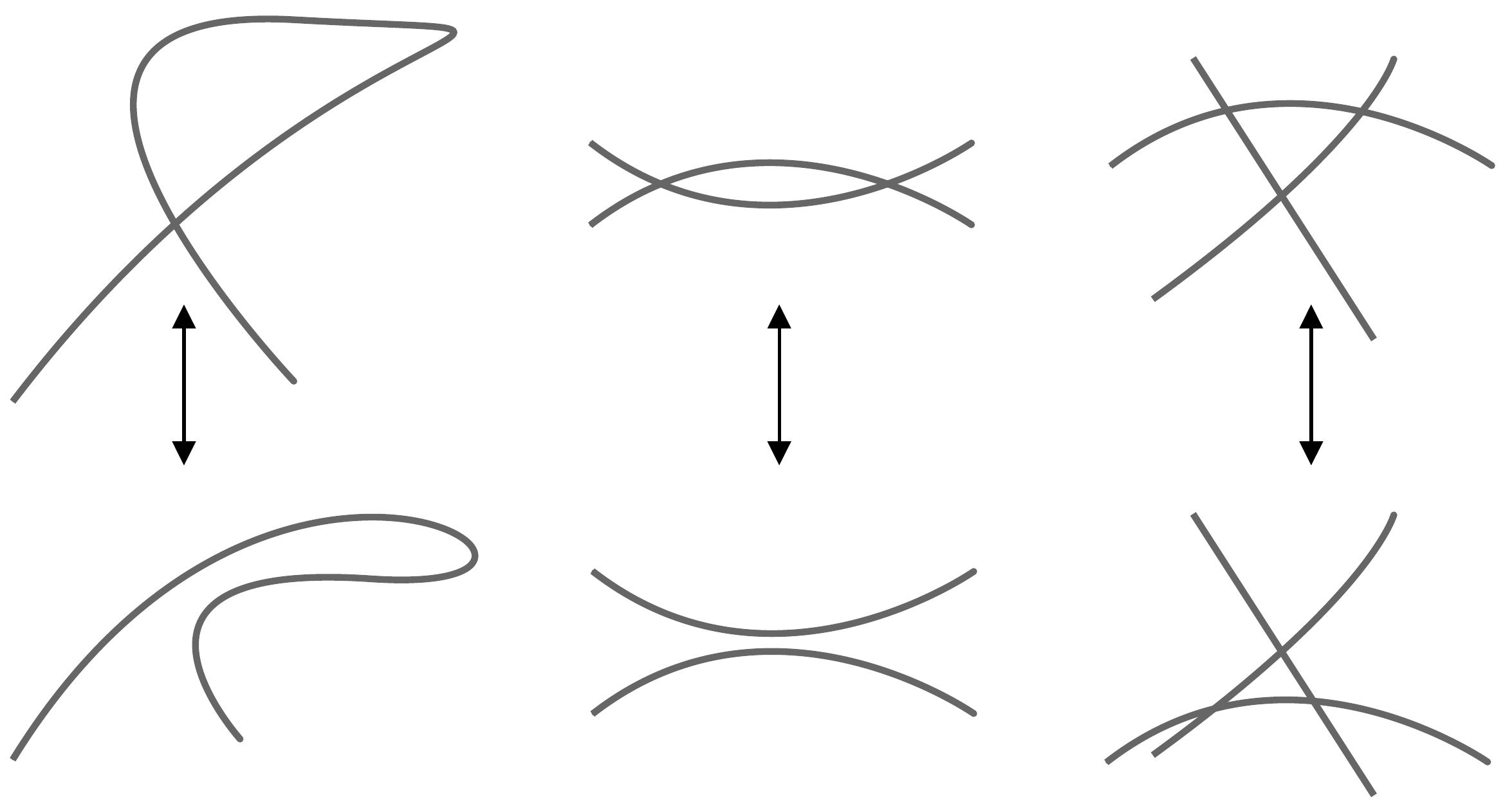}
    \caption{The three types of visual events of a space curve. }
    \label{fig:reidemeister}
\end{figure}

From the algebraic perspective, visual events have been studied by Petitjean \cite{petitjean1996complexity} and Kohn, Sturmfels, and Trager \cite{kohn2018changing}.
The \emph{visual event surface} $\mathcal{V}(X)$ of the curve or surface $X$ is the Zariski closure in $\PP^3$ of the set of camera centers where a visual event occurs.
For a curve $X$ in $\PP^3$, the visual event surface typically has three irreducible components, each representing a different type of visual event.
On the image curve $Y_C(X)$, these types correspond to the three Reidemeister moves from knot theory \cite{reidemeister2013knotentheorie}; see Figure~\ref{fig:reidemeister}

If $X$ is a surface, its visual event surface $\mathcal{V}(X)$ typically has five irreducible components that arise from six different types of visual events (here two types of events cannot be distinguished algebraically, i.e. they lie on the same component of $\mathcal{V}(X)$, but they represent two distinct behaviors over the real numbers).
These six types of events correspond to the non-generic singularities from catastrophe theory \cite{thom1972structural,arnold1984catastrophe}.

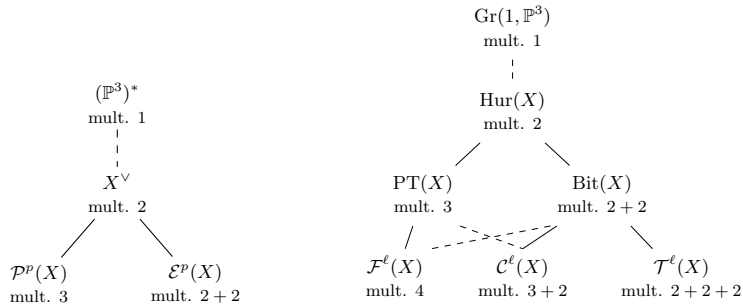
\begin{figure}[htbp]
\small 
  \centering 
  \begin{tikzpicture}[sibling distance=8em, scale=0.8, 
    every node/.style = {%shape=rectangle, rounded corners, draw, 
    align=center, scale=0.8}]] 
    \node {$(\PP^3)^*$\\{\footnotesize mult. $1$}}
    child[dashed] {node[solid] {$X^\vee$\\ \footnotesize mult. $2$}
    child[solid] { node {$\mathcal P^p(X)$\\ \footnotesize mult. $3$} }
    child[solid] { node {$\mathcal E^p(X)$\\ \footnotesize mult. $2+2$} }};
    \end{tikzpicture}
    \qquad \qquad
    \scalebox{0.8}{
    \begin{forest}
    [{${\rm Gr}(1,\PP^3)$ \\ \footnotesize mult. $1$}, align=center, base=bottom
      [{${\rm Hur}(X)$\\ \footnotesize mult. $2$}, align=center, base=bottom, edge=dashed
        [{${\rm PT}(X)$\\\footnotesize mult. $3$}, align =center, base = bottom, name = PTnode
          [{$\mathcal F^\ell(X)$\\\footnotesize mult. $4$}, align=center, base=bottom, name = Fnode]
          [, no edge ]
        ] 
        [{${\rm Bit}(X)$\\\footnotesize mult. $2+2$}, align=center, base=bottom, name = BITnode
          [{$\mathcal C^\ell(X)$\\\footnotesize mult. $3+2$}, align=center, base=bottom, name = Cnode, child anchor = north]
          [, no edge ]
          [{${\mathcal T}^\ell(X)$\\\footnotesize mult. $2+2+2$}, align=center, base=bottom]
        ]
      ]
    ]
    \draw (PTnode.south east)[dashed] -- (Cnode.north);
    \draw (BITnode.south west)[dashed] -- (Fnode.north east);
    \end{forest} }
  \caption{Loci of planes and lines that meet a surface $X$ with assigned multiplicities.
  For instance, ``mult. $2+2+2$'' on the right-hand side refers to the locus of tritangent lines.
  A dashed (resp. solid) edge means that the lower variety is contained in (the singular locus of) the upper one \cite{kohn2018changing}.}
  \label{fig:surface-graph}
\end{figure}

In \cite{kohn2018changing}, the components of the visual event surface $\mathcal{V}(X)$ of a curve or surface $X$ are alternatively characterized via the iterated singular loci of the \emph{higher associated hypersurfaces} \cite{gelfand2008discriminants} / \emph{coisotropic hypersurfaces} \cite{kohn2021isotropic} of $X$. 
For a surface $X$, the coisotropic hypersurfaces are the \emph{dual surface} $X^\vee$ in $(\PP^3)^\ast$ 
and the \emph{Hurwitz threefold} $\mathrm{Hur}(X)$ in $\Gr(1, \PP^3)$,
i.e., the Zariski closure of the set of planes (respectively lines) tangent to $X$.
Figure \ref{fig:surface-graph} shows the iterated singular loci of these hypersurfaces till the level of curves. 
The five components of $\mathcal{V}(X)$ are encoded by the curves in the last row of Figure \ref{fig:surface-graph}:
two components are the dual surfaces to the curves on the left. The remaining three components are ruled by the lines on the curves shown in the right diagram.

\begin{figure}[htb]
    \centering
    \begin{subfigure}[b]{0.31\textwidth}
         \centering
         \includegraphics[width=\textwidth]{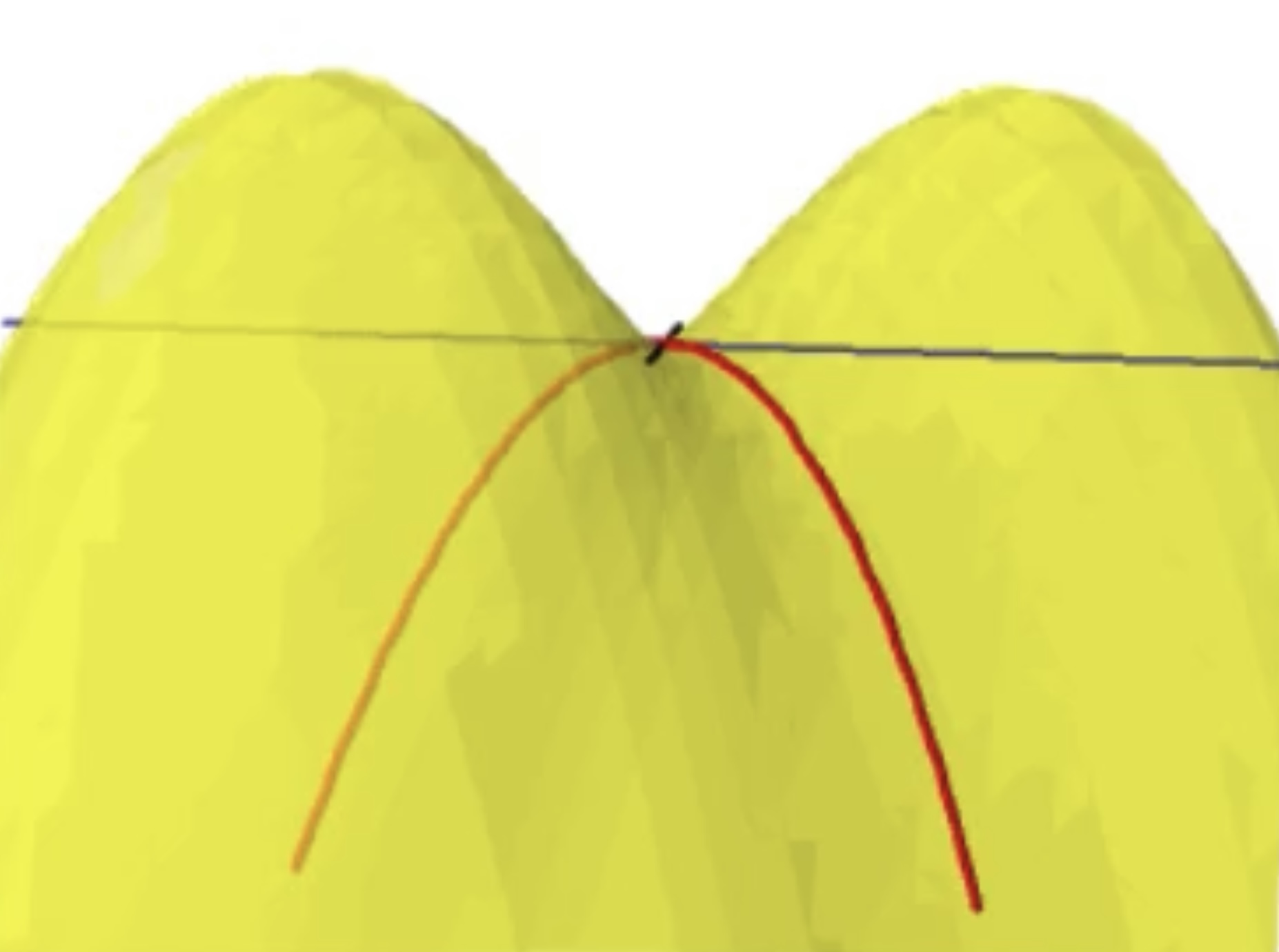}
         \caption{Smooth image contour}
         \label{fig:swallowtail1}
     \end{subfigure}
          \begin{subfigure}[b]{0.3\textwidth}
         \centering
         \includegraphics[width=\textwidth]{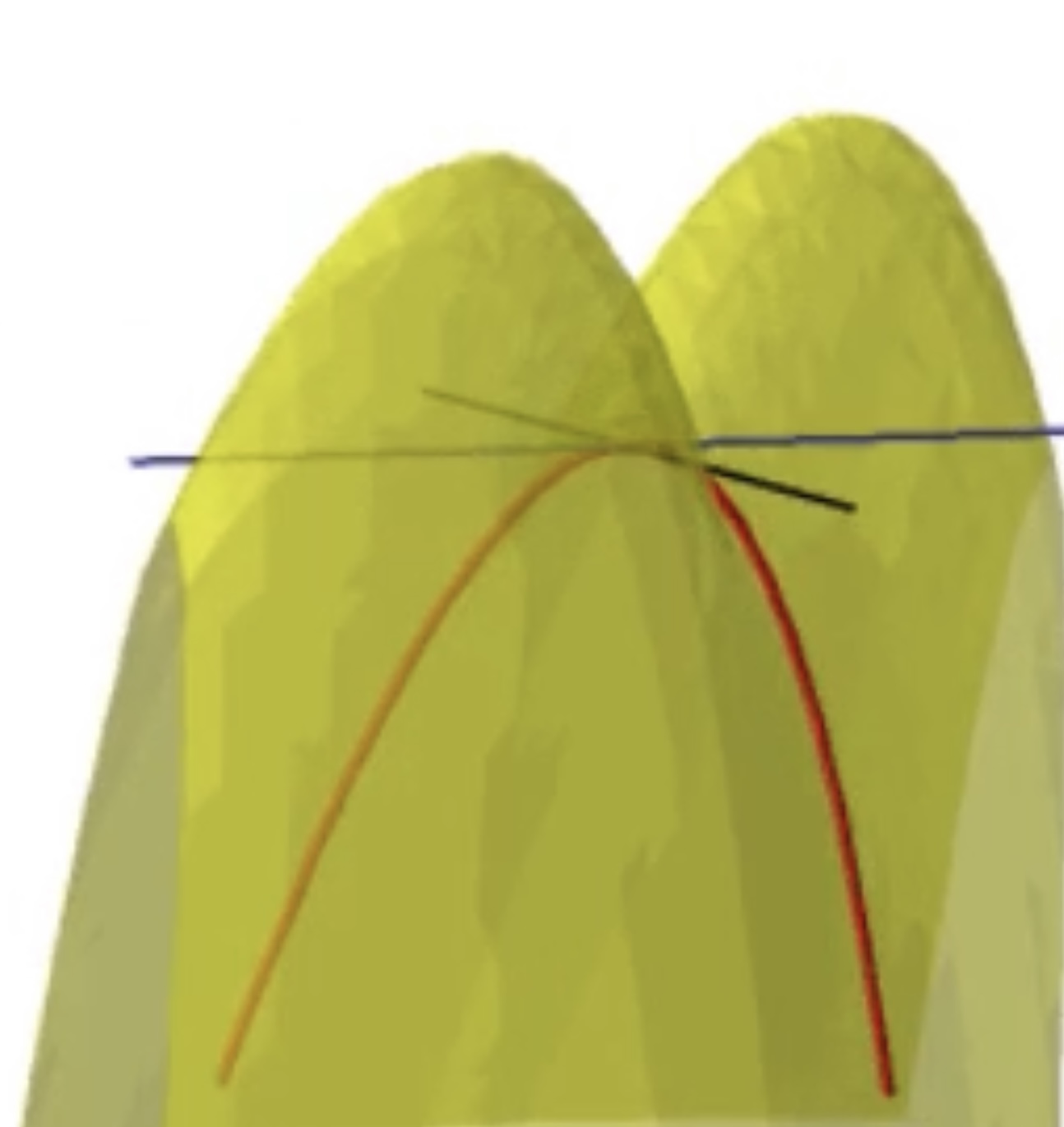}
         \caption{Singular image contour}
         \label{fig:swallowtail2}
     \end{subfigure}
          \begin{subfigure}[b]{0.31\textwidth}
         \centering
         \includegraphics[width=\textwidth]{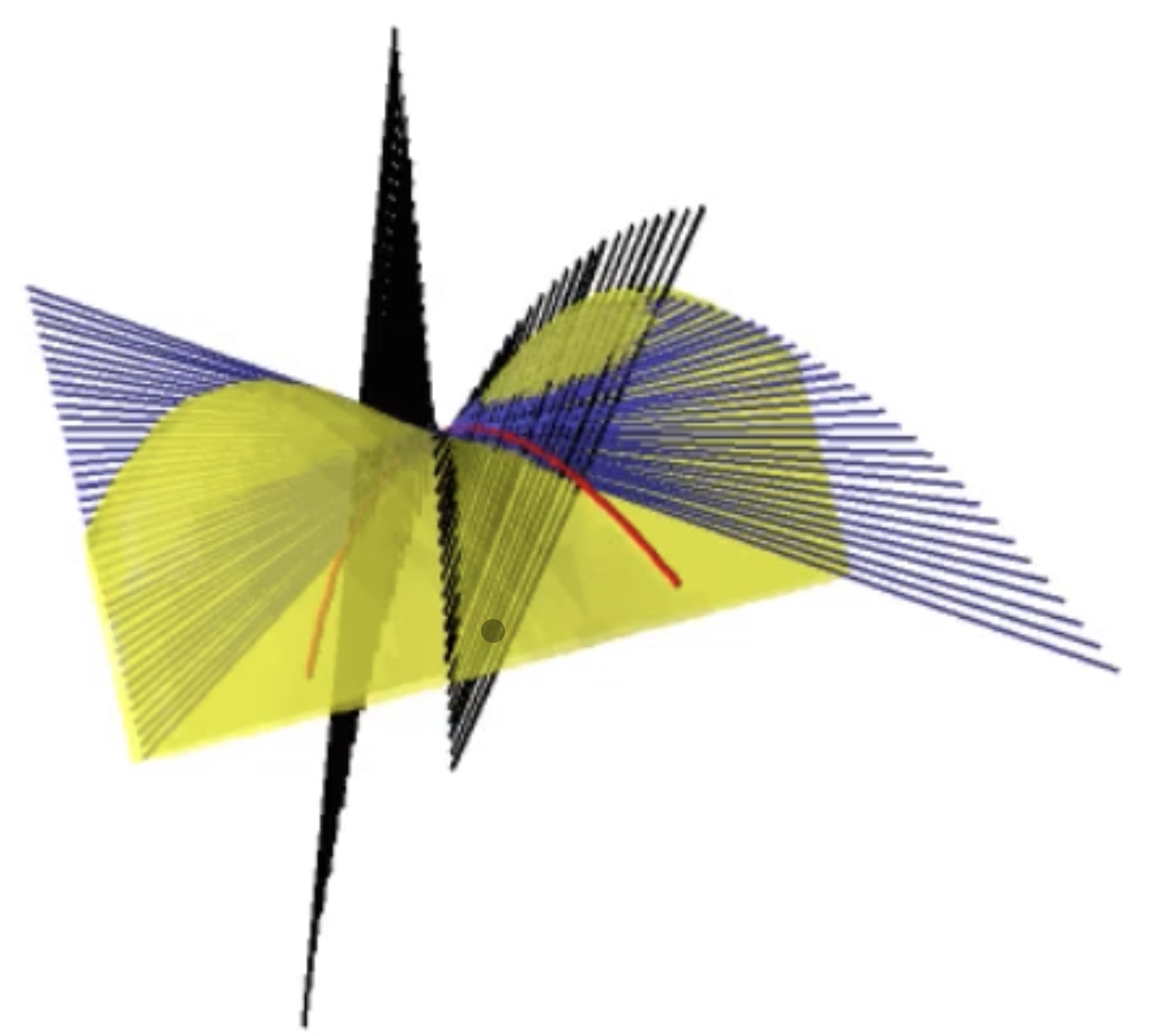}
         \caption{Flecnodal surface}
         \label{fig:flecnodal}
     \end{subfigure}
    \caption{Quartic surface parametrized by $(s,t,s^2+3st+t^4)$ (in yellow) with flecnodal lines that meet the surface with multiplicity $4$ at one of the red points.
    A swallowtail event happens when the $2$ cusps and the node in Figure \ref{fig:swallowtail2} come  together to a higher order singularity.}
    \label{fig:swallowtail}
\end{figure}

\begin{example}
    Consider the curve $\mathcal{F}^\ell(X) \subset \mathrm{Gr}(1, \PP^3)$ associated with a surface $X \subset \PP^3$ of degree at least four.
    The curve $\mathcal{F}^\ell(X)$ is the Zariski closure of the set of all lines that intersect the surface $X$ at some point with multiplicity four. 
    Such lines are called \emph{flecnodal}.
    The union of all flecnodal lines is the \emph{flecnodal surface} of $X$; it is one of the irreducible components of the visual event surface $\mathcal{V}(X)$. The visual event corresponding to the flecnodal surface is called \emph{swallowtail event}.

    For an explicit example, let us consider the quartic surface parametrized by $(s,t) \mapsto (s,t,s^2+3st+t^4)$. 
    The surface is shown in yellow in Figure~\ref{fig:swallowtail}. 
    In Figure \ref{fig:swallowtail1}, the surface's image contour is smooth.  
    In Figure \ref{fig:swallowtail2}, the image contour has two cusps and a node, as illustrated in the torus example in Figure \ref{fig:contour}.
    The transition occurs when the camera center is located on a flecnodal line. 
    On a general surface, there is a one-dimensional family of points such that one of the tangent lines at that point is flecnodal. The family makes up the curve shown in red in Figure \ref{fig:swallowtail}.
    For the particular surface there, every red point has two flecnodal lines.
    The lines are shown for a particular point in Figures \ref{fig:swallowtail1} and \ref{fig:swallowtail2}, and for many points in Figure \ref{fig:flecnodal}.
\end{example}

Similarly, if $X$ is a curve its coisotropic hypersurfaces are the dual surface $X^\vee$ in $(\PP^3)^\ast$ and the \emph{Chow threefold} in $\Gr(1, \PP^3)$, that is, the set of lines meeting $X$. Their (iterated) singular curves give rise to the three components of $\mathcal{V}(X)$ analogously.

The degrees of the components of $\mathcal{V}(X)$ for a general curve or surface $X$ of fixed degree (and genus in the curve case) are classically known; see \cite[Theorems 3.1 and 4.1]{kohn2018changing}.
Petitjean provided degree formulas for the five components of $\mathcal{V}(X)$ in the surface case \cite{petitjean1996complexity}, 
but they appeared in fact already in paragraphs 597, 598, 599, 608 and 613 of Salmon’s 1882 book \cite{salmon1912treatise}.

\subsection{Hurwitz hypersurfaces and condition numbers}
\label{ssec:hurwitz}

The minimal problem \eqref{eq:solveE} of computing essential matrices consistent with five point pairs admits the following geometric interpretation: 
We are computing the intersection of the essential variety $\mathcal{E} \subset \PP(\RR^{3 \times 3})$ with a data-dependent linear subspace $L \subset \PP(\RR^{3 \times 3})$ of complementary dimension, namely $L = \{ x_i y_i^{\top} : i \in \mathcal{I} \}^{\perp}$. 
In numerical algebraic geometry language, we are computing a \textit{witness set} of $\mathcal{E}$ \cite[Chapter 13.3]{sommese2005numerical}.
Several, though not all, minimal problems in computer vision usefully admit such a interpretation.  
For instance, Kileel \cite{kileel2017minimal} describes all $66$ minimal problems for three calibrated pinhole cameras that arise from slicing the calibrated trifocal variety $\mathcal{T}_{\mathrm{cal}}$ with linear spaces of complementary dimension.

In the applied algebra community, this description has motivated theoretical works on intersecting a fixed variety $Z \subset \PP^n$ with varying linear subspaces of complementary dimension.
This includes the detailed study of the \emph{Hurwitz hypersurface} $\mathrm{Hur}(Z)$ by Sturmfels \cite{sturmfels2017hurwitz}.
Writing $c$ for the codimension of $Z$ in $\PP^n$ and $d$ for its degree, $\mathrm{Hur}(Z)$ is the Zariski closure in $\Gr(c, \PP^n)$ of the set of all $c$-dimensional projective subspaces $L$ of $\PP^n$ such that the intersection $Z \cap L$ does not consist of $d$ reduced points.
The Hurwitz hypersurface is intimately linked with the \emph{condition number} \cite{burgisser2013condition} of the algebraic function $L \mapsto Z \cap L$.
That condition number measures how much the intersection $Z \cap L$ changes when $L$ gets perturbed, which is crucial for the understanding of how much numerical computations such as homotopy continuation are affected by errors.
B\"urgisser  showed that the set of \emph{ill-posed} $L$, i.e. those with infinite condition, is essentially the Hurwitz hypersurface $\mathrm{Hur}(Z)$ and that (almost) all $L$ with large condition are contained in a small tube around $\mathrm{Hur}(Z)$ \cite[Corollary 1.10]{burgisser2017condition}.
B\"urgisser's work follows the general idea that condition numbers are given by the inverse distance to ill-posedness from numerical analysis \cite{demmel1987geometry,demmel1987condition}.

Recently there has been work on condition numbers of minimal problems in the more practical setting, where the dependence is in terms of the perturbation of the image data (instead of perturbation of $L$ in the Grassmannian) by Fan, Kileel and Kimia \cite{fan2021instability, fan2023condition}.  
Their work drew distinctions between criticality (Section~\ref{ssec:critical} above) and ill-posedness (when the condition number is infinite), see  \cite[Section~5]{fan2023condition}.   
For the $5$- and $7$-point problems, the authors derived condition number formulas in terms of certain Jacobian matrices, and they characterized ill-posed world scenes and ill-posed image point pairs geometrically.
Furthermore, real-data experiments showed that poor conditioning can actually plague some computer vision data sets.


\bibliographystyle{amsplain}
\bibliography{literature}

\end{document}